\documentclass[11pt]{amsart}

\usepackage{epigamath}

\usepackage[english]{babel}

\numberwithin{equation}{section}

\usepackage[shortlabels,inline]{enumitem}

\usepackage{babel}
\usepackage{amsbsy}
\usepackage{amstext}
\usepackage{amsthm}
\usepackage{amssymb}
\usepackage{setspace}

\theoremstyle{plain}
\newtheorem{thm}{Theorem}[section]
\newtheorem{prop}[thm]{Proposition}
\newtheorem{lem}[thm]{Lemma}

\theoremstyle{definition}
\newtheorem{defn}[thm]{Definition}

\theoremstyle{remark}
\newtheorem*{claim*}{Claim}
\newtheorem{claim}{Claim}

\setlist[enumerate,1]{label={\rm(\arabic*)}, ref={\rm\arabic*}} 

\DeclareMathOperator{\diam}{diam}
\DeclareMathOperator{\Gal}{Gal}
\DeclareMathOperator{\PrePer}{PrePer}

\newcommand{\A}{\mathbb A}
\newcommand{\C}{\mathbb C}
\newcommand{\PP}{\mathbb P}
\newcommand{\Q}{\mathbb Q}
\newcommand{\R}{\mathbb R}

\newcommand{\an}{\mathrm{an}}
\newcommand{\Diag}{\mathrm{Diag}}
\newcommand{\bounded}{\mathrm{bounded}}
\newcommand{\close}{\mathrm{close}}
\newcommand{\help}{\mathrm{help}}

\newcommand{\supth}[1]{\ensuremath{#1^{\mathrm{th}}}}

\renewcommand{\bar}{\overline}

\EpigaVolumeYear{10}{2026} \EpigaArticleNr{15} \ReceivedOn{February 14, 2023}
\InFinalFormOn{December 17, 2025}
\AcceptedOn{March 2, 2026}

\title{Uniform unlikely intersections for unicritical polynomials}
\titlemark{Uniform unlikely intersections for unicritical polynomials}

\author{Hang Fu}
\address{Department of Mathematics and Computer Science, University of Basel,
Spiegelgasse 1, 4051 Basel, Switzerland}
\email{drfuhang@gmail.com}

\authormark{H.~Fu}

\AbstractInEnglish{Fix $d\geq2$, and consider the family
  $f_{t}(z)=z^{d}+t$ of polynomials parameterized by $t\in\C$. In this
  article, we will show that there exists a constant $C(d)$ such that
  for any $a,b\in\C$ with $a^{d}\neq b^{d}$, the number of $t\in\C$
  such that $a$ and $b$ are both preperiodic for $f_{t}$ is at most
  $C(d)$.}


\begin{document}



\maketitle

\begin{prelims}

\DisplayAbstractInEnglish

\bigskip

\DisplayKeyWords

\medskip

\DisplayMSCclass

\end{prelims}


\newpage

\setcounter{tocdepth}{1}

\tableofcontents


\section{Introduction}\label{sec1}

Fix an integer $d\geq2$, and consider the family 
\[
f_{t}(z)=z^{d}+t
\]
of polynomials parameterized by $t\in\C$. A
point $a\in\C$ is said to be \textit{preperiodic} for $f_{t}$
if its forward orbit $O_{f_{t}}(a)=\{f_{t}^{n}(a):n\geq1\}$ is finite,
where $f_{t}^{n}(a)$ is the $\supth{n}$ iterate of $a$ under $f_{t}$.

Let $\PrePer(f_{t})$ be the set of all preperiodic points of
$f_{t}$. By \cite[Section 4]{Be}, \cite[Theorem 1.2 and Corollary 1.3]{BD1},
and \cite[Theorem 1.3]{YZ}, the intersection $\PrePer(f_{t_{1}})\cap\PrePer(f_{t_{2}})$
is finite if and only if $t_{1}\neq t_{2}$. Furthermore, when $d=2$,
DeMarco, Krieger, and Ye \cite[Theorem 1.1]{DKY2} proved that $\PrePer(f_{t_{1}})\cap\PrePer(f_{t_{2}})$
is in fact uniformly bounded for any $t_{1}\neq t_{2}$. This result
is an analog of \cite[Theorem 1.4]{DKY1}, which provides a partial
solution to the effective finiteness conjecture for elliptic curves
proposed in \cite{BFT} and also implies a uniform Manin--Mumford
bound under suitable conditions; see \cite[Theorem 1.1]{DKY1}.

On the other hand, for any $a,b\in\C$, we define
\[
S_{a,b}=\{t\in\C:a\text{ and }b\text{ are both preperiodic for }f_{t}\}.
\]
Note that $S_{a,b}$ depends on $d$ implicitly. Zannier \cite[Section 3.4.7]{Za}
asked whether $S_{0,1}$ is finite when $d=2$. Baker and DeMarco
\cite[Theorem 1.1]{BD1} answered this question affirmatively by showing
that $S_{a,b}$ is finite if and only if $a^{d}\neq b^{d}$. This
result is motivated by a theorem of Masser and Zannier. In \cite{MZ1,MZ2,MZ3}, they showed that for any $a\neq b\in\C\backslash\{0,1\}$,
there exist only finitely many $t\in\C$ such that
\[
\left(a,\sqrt{a(a-1)(a-t)}\right)\quad\text{and}\quad\left(b,\sqrt{b(b-1)(b-t)}\right)
\]
are both torsion on the elliptic curve $y^{2}=x(x-1)(x-t)$. For further
developments on this problem, see also \cite{BD2,GHT1,GHT2}.

However, the proof of \cite[Theorem 1.1]{BD1} is not effective, so no
explicit upper bound for $|S_{a,b}|$ is given. In order to find more
information on $S_{a,b}$, Fili \cite{Fi} studied the case $S_{0,1}$
when $d=2$ and, more importantly, made a key observation, see
\cite[Theorem 1]{Fi}, which turns out to be crucial for \cite{DKY1}
and \cite{DKY2}. In light of \cite[Theorem 1.4]{DKY1} and
\cite[Theorem 1.1]{DKY2}, one is led to ask whether $S_{a,b}$ is also
uniformly bounded for any $a,b\in\C$ with $a^{d}\neq b^{d}$. In this
article, we are able to show that the answer is yes.

\begin{thm}\label{thm1.1}
For any integer $d\geq2$, there exists a constant
$C(d)$ such that $|S_{a,b}|\leq C(d)$ for any $a,b\in\C$
with $a^{d}\neq b^{d}$.
\end{thm}

Our results and proofs are inspired by \cite{DKY2}, and our techniques
and strategies also come from their ideas. By a standard specialization
argument, it suffices to prove Theorem~\ref{thm1.1} for $a,b\in\bar{\Q}$.
Let $K$ be a number field such that $a,b\in K$, and let $M_{K}$
be the set of places. For each $v\in M_{K}$, we work with the dynamics
of $f_{t}$ on the Berkovich affine line $\A_{v}^{1,\an}$.
In Section~\ref{sec2}, we will introduce
\begin{itemize}
\item the generalized Mandelbrot sets $M_{a,v}$ and $M_{b,v}$,
\item the equilibrium measures $\mu_{a,v}$ and $\mu_{b,v}$,
\item the Green's functions $g_{a,v}$ and $g_{b,v}$.
\end{itemize}
Then the Arakelov--Zhang pairing of $\boldsymbol{\mu}_{a}=\{\mu_{a,v}\}_{v\in M_{K}}$
and $\boldsymbol{\mu}_{b}=\{\mu_{b,v}\}_{v\in M_{K}}$ is given by
\[
\left\langle \boldsymbol{\mu}_{a},\boldsymbol{\mu}_{b}\right\rangle =\sum_{v\in M_{K}}\frac{[K_{v}:\Q_{v}]}{[K:\Q]}\int_{\A_{v}^{1,\an}}g_{a,v}d\mu_{b,v}.
\]
The value of $\langle \boldsymbol{\mu}_{a},\boldsymbol{\mu}_{b}\rangle $
depends on $a$ and $b$ only and is independent of the choice of
$K$. By \cite[Propositions 2.6 and 4.5]{FRL} and \cite[Theorems 1.1 and 3.4]{BD1}, we have 
$\langle \boldsymbol{\mu}_{a},\boldsymbol{\mu}_{b}\rangle =\langle \boldsymbol{\mu}_{b},\boldsymbol{\mu}_{a}\rangle \geq0$
and
\[
\left\langle \boldsymbol{\mu}_{a},\boldsymbol{\mu}_{b}\right\rangle =0\Longleftrightarrow\boldsymbol{\mu}_{a}=\boldsymbol{\mu}_{b}\Longleftrightarrow a^{d}=b^{d}\Longleftrightarrow|S_{a,b}|=\infty.
\]
Given this equivalence relation, it is not surprising that the value
of $\langle \boldsymbol{\mu}_{a},\boldsymbol{\mu}_{b}\rangle $
encodes some information about $S_{a,b}$. As in \cite{DKY1, DKY2}, the main task of this article is to estimate the upper
and lower bounds for $\langle \boldsymbol{\mu}_{a},\boldsymbol{\mu}_{b}\rangle $.
More precisely, we prove the following counterparts of \cite[Theorems 1.5, 1.6, and~1.7]{DKY1}
and \cite[Theorems 1.6, 1.7, and 1.9]{DKY2}.

\begin{thm}\label{thm1.2}
Let $a,b\in\bar{\Q}$ such that $a^{d}\neq b^{d}$
and $|S_{a,b}|>0$. For any $0<\varepsilon<4d$, we have
\[
\left\langle \boldsymbol{\mu}_{a},\boldsymbol{\mu}_{b}\right\rangle \leq\left(\varepsilon+\frac{8d/\varepsilon-2}{|S_{a,b}|}\right)(h(a,b)+5),
\]
where $h$ is the logarithmic Weil height on $\A^{2}(\bar{\Q})$.
\end{thm}

\begin{thm}\label{thm1.3}
  For any $a,b\in\bar{\Q}$ such that $a^{d}\neq b^{d}$,
we have
\[
\left\langle \boldsymbol{\mu}_{a},\boldsymbol{\mu}_{b}\right\rangle \geq\frac{1}{12d^{2}}h(a,b)-1,
\]
where $h$ is the logarithmic Weil height on $\A^{2}(\bar{\Q})$.
\end{thm}

\begin{thm} \label{thm1.4}
  There exists a constant $\delta(d)>0$ such that $\langle \boldsymbol{\mu}_{a},\boldsymbol{\mu}_{b}\rangle \geq\delta(d)$
for any $a,b\in\bar{\Q}$ with $a^{d}\neq b^{d}$.
\end{thm}

Note that Theorem~\ref{thm1.3} is trivial when $h(a,b)$ is small.
We need Theorem~\ref{thm1.4} to deal with this case.

The main differences between this article and \cite{DKY2} are:
\begin{enumerate*}\item
Our Theorem~\ref{thm1.1} is valid for any $d\geq2$, while
\cite[Theorem 1.1]{DKY2} focuses on $d=2$. \item We work with the
generalized Mandelbrot sets instead of the Julia sets in
\cite{DKY2}. \item When we estimate the lower bounds for $\left\langle
\mu_{a,v},\mu_{b,v}\right\rangle _{v}$ at the non-Archimedean places
$v\in M_{K}^{0}$, our computations in Section~\ref{sec4.4} are more
simplified than those in \cite[Sections 5 and~6]{DKY2}.  This
simplification helps us to work with all $d\geq2$ at the same
time.\end{enumerate*}

The plan of this article is as follows: In Section~\ref{sec2}, we fix
the notation and review the tools we will need. In Section~\ref{sec3},
we first estimate the upper bounds for $\langle
\boldsymbol{\mu}_{a},\boldsymbol{\mu}_{b}\rangle $ locally, and then
we combine the local estimates to give the proof of
Theorem~\ref{thm1.2}. The structure of Section~\ref{sec4} is similar,
but this time we estimate the lower bounds for $\langle
\boldsymbol{\mu}_{a},\boldsymbol{\mu}_{b}\rangle $ and give the proofs
of Theorems~\ref{thm1.3} and~\ref{thm1.4}. In Section~\ref{sec5}, we
show that Theorem~\ref{thm1.1} follows from
Theorems~\ref{thm1.2},~\ref{thm1.3}, and~\ref{thm1.4}.

\subsection*{Acknowledgments} The author thanks Liang-Chung Hsia for
many insightful discussions, for his careful reading of earlier drafts,
and for his valuable comments and suggestions. The author thanks the
anonymous referees for their helpful feedback. The author thanks National
Taiwan Normal University for the hospitality.

\section{Notation and preliminaries}\label{sec2}

The main references for this section are \cite{BR2,BD1,FRL,Fi}.

Given a number field $K$, let $M_{K}$ be the set of places, let
$M_{K}^{\infty}$ be the set of Archimedean places, and let $M_{K}^{0}$
be the set of non-Archimedean places. We normalize the absolute values
$|\cdot|_{v}$ on $K$ such that they extend the standard absolute
values on $\Q$. For any $a\in K^{\times}$, we have the product
formula
\[
\prod_{v\in M_{K}}|a|_{v}^{n_{v}}=1,\quad\text{where }n_{v}=\frac{[K_{v}:\Q_{v}]}{[K:\Q]}.
\]

\subsection{Berkovich spaces}\label{sec2.1}

For each $v\in M_{K}$, let $K_{v}$ be the completion of $K$ at $v$,
let $\bar{K}_{v}$ be an algebraic closure of $K_{v}$, and let $\C_{v}$
be the completion of $\bar{K}_{v}$. The \textit{Berkovich affine line}
$\A_{v}^{1,\an}$ is a locally compact, Hausdorff, path-connected space
containing $\C_{v}$ as a dense subspace. As a topological space,
$\A_{v}^{1,\an}$ is the set of all multiplicative seminorms
$[\cdot]_{x}\colon\C_{v}[T]\to\R$ on the polynomial ring $\C_{v}[T]$
that extend the absolute value $|\cdot|_{v}$ on $\C_{v}$, endowed with
the weakest topology for which $x\mapsto[f]_{x}$ is continuous for any
$f\in\C_{v}[T]$.  The \textit{Berkovich projective line}
$\PP_{v}^{1,\an}$ can be identified with the one-point
compactification of $\A_{v}^{1,\an}$.

If $v\in M_{K}^{\infty}$, then by the Gelfand--Mazur theorem,
$\A_{v}^{1,\an}$ is homeomorphic to $\C_{v}=\C$.
If $v\in M_{K}^{0}$, then by Berkovich\textquoteright s classification
theorem, each $x\in\A_{v}^{1,\an}$ corresponds to a
decreasing nested sequence $\{\overline{D}(a_{n},r_{n})\}_{n=1}^{\infty}$
of closed disks on $\C_{v}$ such that
\[
[f]_{x}=\lim_{n\to\infty}\sup_{z\in\overline{D}(a_{n},r_{n})}|f(z)|_{v}.
\]
Based on the nature of $D=\cap_{n=1}^{\infty}\overline{D}(a_{n},r_{n})$,
the points of $\A_{v}^{1,\an}$ can be categorized into
four types: \begin{enumerate*}[label=(\Roman*)]\item $D$ is a point of $\C_{v}$, \item $D$ is
a closed disk with radius in $|\C_{v}^{\times}|$, \item $D$
is a closed disk with radius not in $|\C_{v}^{\times}|$,
and \item $D$ is the empty set.\end{enumerate*}

For any $a\in\C_{v}$ and any $r>0$, we define $\mathcal{D}(a,r)$ to be
the set of points corresponding to
$\{\overline{D}(a_{n},r_{n})\}_{n=1}^{\infty}$ with
$\overline{D}(a_{n},r_{n})\subseteq\overline{D}(a,r)$ and define
$\zeta_{a,r}$ to be the point corresponding to $\overline{D}(a,r)$.

\subsection{Potential theory}\label{sec2.2}

When $v\in M_{K}^{0}$, we introduce the \textit{Hsia kernel}
\[
\delta_{v}(x,y)=\limsup_{z,w\in\C_{v},z\to x,w\to y}|z-w|_{v},
\]
which extends the distance function $|x-y|_{v}$ on $\C_{v}$ to the
entire $\A_{v}^{1,\an}$. When $v\in M_{K}^{\infty}$, we also write
$\delta_{v}(x,y)=|x-y|_{v}$ to unify the notation further on.

Fix a $v\in M_{K}$, and let $E$ be a compact subset of $\A_{v}^{1,\an}$.
The \textit{logarithmic capacity} $\gamma_{v}(E)$ of $E$ is given
by
\[
-\log\gamma_{v}(E)=\inf_{\mu}\iint_{E\times E}-\log\delta_{v}(z,w)d\mu(z)d\mu(w),
\]
where the infimum is taken over all probability measures $\mu$ supported
on $E$. If $\gamma_{v}(E)>0$, then there exists a unique probability
measure $\mu_{E}$, called the \textit{equilibrium measure} of $E$,
such that the infimum is achieved. The \textit{Green's function} of
$E$ is defined by
\[
g_{E}(z)=-\log\gamma_{v}(E)+\int_{E}\log\delta_{v}(z,w)d\mu_{E}(w),
\]
which is a non-negative real-valued function on $\A_{v}^{1,\an}$.

\subsection{Generalized Mandelbrot sets}\label{sec2.3}

Now we go back to the dynamics of $f_{t}(z)=z^{d}+t$. For any $v\in M_{K}$
and any $a\in\C_{v}$, we define the \textit{generalized Mandelbrot
set} by
\[
M_{a,v}=\left\{t\in\A_{v}^{1,\an}:\sup_{n}\left[f_{T}^{n}(a)\right]_{t}<\infty\right\},
\]
where $f_{T}^{n}(a)$ is considered as an element of the polynomial
ring $\C_{v}[T]$. Note that if $t\in\C_{v}$, then
$[f_{T}^{n}(a)]_{t}$ is simply $|f_{t}^{n}(a)|_{v}$. It is known
that $M_{a,v}$ is compact. We write $\mu_{a,v}$ and $g_{a,v}$ for
the equilibrium measure and the Green's function of $M_{a,v}$. The
following properties are collected from \cite[Section 3]{BD1}.

\begin{thm}
\label{thm2.1}For any $v\in M_{K}$ and any $a\in\C_{v}$,
we have the following:
\begin{enumerate}
\item The logarithmic capacity $\gamma_{v}(M_{a,v})$ equals $1$.
\item The Green's function of\, $M_{a,v}$ is given by
\[
g_{a,v}(t)=\lim_{n\to\infty}\frac{1}{d^{n}}\log^{+}\left[f_{T}^{n+1}(a)\right]_{t},
\]
where $\log^{+}z=\log\max\{z,1\}$ for any $z\in\R$.
\item The function $g_{a,v}(t)$ is continuous on $\A_{v}^{1,\an}$.
\item The function $g_{a,v}(t)$ is harmonic on $\A_{v}^{1,\an}\backslash M_{a,v}$.
\item We have  $g_{a,v}(t)=0$ if and only if $t\in M_{a,v}$.
\end{enumerate}
\end{thm}

\subsection{Arakelov--Zhang pairing}\label{sec2.4}

Like \cite[Theorem 1.7]{DKY1} and \cite[Theorem 1.9]{DKY2}, our Theorem
\ref{thm1.2} also builds on the quantitative equidistribution results
of \cite{FRL,Fi}.

Let $a\in\C_{v}$ and $r>0$. If $v\in M_{K}^{\infty}$, we
define $m_{a,r,v}$ to be the normalized Haar measure on the circle
$\partial D(a,r)$. If $v\in M_{K}^{0}$, we define $m_{a,r,v}$ to
be the Dirac measure on $\zeta_{a,r}$.

\begin{defn}[\textit{cf.} {\cite[D\'efinition 1.1]{FRL}}]\label{def2.2}
We call $\boldsymbol{\mu}=(\mu_{v})_{v\in M_{K}}$
an \textit{adelic measure} if
\begin{enumerate}
\item $\mu_{v}$ is a probability measure on $\PP_{v}^{1,\an}$
for any $v\in M_{K}$,
\item $\mu_{v}=m_{0,1,v}$ for all but finitely many $v\in M_{K}$,
\item for any $v\in M_{K}$, we have $\mu_{v}-m_{0,1,v}=\Delta u_{v}$ for some
continuous function $u_{v}$ on $\PP_{v}^{1,\an}$, where
$\Delta$ is the Laplacian on $\PP_{v}^{1,\an}$.
\end{enumerate}
\end{defn}

Following \cite[Sections 2.4 and 4.4]{FRL}, for each $v\in M_{K}$,
we define the \textit{mutual energy} of two signed measures $\mu_{1,v}$
and $\mu_{2,v}$ on $\PP_{v}^{1,\an}$ by
\[
\left(\mu_{1,v},\mu_{2,v}\right)_{v}=\iint_{\A_{v}^{1,\an}\times\A_{v}^{1,\an}\backslash\Diag_{v}}-\log\delta_{v}(z,w)d\mu_{1,v}(z)d\mu_{2,v}(w),
\]
where $\Diag_{v}$ is the diagonal on $\C_{v}\times\C_{v}$.
Suppose $\boldsymbol{\mu}_{1}$ and $\boldsymbol{\mu}_{2}$ are adelic
measures; then we define their \textit{$v$-adic Arakelov--Zhang
pairing} by
\[
\left\langle \mu_{1,v},\mu_{2,v}\right\rangle _{v}=\frac{1}{2}\left(\mu_{1,v}-\mu_{2,v},\mu_{1,v}-\mu_{2,v}\right)_{v} 
\]
and define their \textit{Arakelov--Zhang pairing} by
\[
\left\langle \boldsymbol{\mu}_{1},\boldsymbol{\mu}_{2}\right\rangle =\sum_{v\in M_{K}}n_{v}\left\langle \mu_{1,v},\mu_{2,v}\right\rangle _{v}.
\]

\begin{thm}[\textit{cf.} {\cite[Theorem 1]{Fi}}]\label{thm2.3}
The square root of the Arakelov--Zhang
pairing $\left\langle \cdot,\cdot\right\rangle ^{1/2}$ gives a metric
on the space of all adelic measures.
\end{thm}

Now we go back to the dynamics of $f_{t}(z)=z^{d}+t$. Given $a,b\in K$,
let $\boldsymbol{\mu}_{a}=\{\mu_{a,v}\}_{v\in M_{K}}$ and $\boldsymbol{\mu}_{b}=\{\mu_{b,v}\}_{v\in M_{K}}$
be the equilibrium measures defined in Section~\ref{sec2.3}. It is
known that they are adelic measures, and their $v$-adic Arakelov--Zhang
pairing can be written as
\[
\left\langle \mu_{a,v},\mu_{b,v}\right\rangle _{v}=\int_{\A_{v}^{1,\an}}g_{a,v}d\mu_{b,v}.
\]

For the Arakelov--Zhang pairing in more general settings,
see \cite{Zh,PST,CL2}.

\section{Upper bounds}\label{sec3}

The purpose of this section is to give the proof of Theorem~\ref{thm1.2}.
In order to do so, we estimate the upper bounds for $g_{a,v}(s)$
when $s$ is close to $M_{a,v}$. We work with $v\in M_{K}^{\infty}$
in Section~\ref{sec3.1} and $v\in M_{K}^{0}$ in Section~\ref{sec3.2}.
In Section~\ref{sec3.3}, we apply the local estimates in Theorem
\ref{thm3.13} to complete the proof of Theorem~\ref{thm1.2}.

\subsection{Archimedean estimates}\label{sec3.1}

In this section, we assume that $K$ is a number field, $a\in K$,
and $v\in M_{K}^{\infty}$. Because $v$ is fixed, we write $|\cdot|$,
$M_{a}$, and $g_{a}$ for $|\cdot|_{v}$, $M_{a,v}$, and $g_{a,v}$.

From \cite[Lemma 3.2]{BD1}, we know that $M_{a}$ is bounded. Their
proof can be modified slightly to give an explicit bound for $M_{a}$.
We begin with a basic distortion result for univalent maps, which
is used in the proof of \cite[Lemma 3.2]{BD1} and also in \cite[Section 3.1]{DKY2}.

\begin{thm}\label{thm3.1}
Let $U_{R}=\{z\in\C:|z|>R\}$. If $\phi\colon U_{R}\to\C$
is analytic and injective, and
\[
\phi(z)=z+\sum_{n=1}^{\infty}\frac{a_{n}}{z^{n}},
\]
then $\phi(U_{R})\supseteq U_{2R}$. In particular, $|\phi(z)|\leq2|z|$
for any $z\in U_{R}$.
\end{thm}

\begin{proof}
The first assertion is \cite[Corollary 3.3]{BH}. If the second assertion
is false, then there exists a $z\in U_{R}$ such that $\phi(z)\in U_{2|z|}\subseteq\phi(U_{|z|})$,
which contradicts the injectivity of $\phi$.
\end{proof}

\begin{prop}\label{prop3.2}
If\,  $|t|>4^{d}\max\{|a|,4\}^{d}$, then
\[
\log|t|-1\leq g_{a}(t)\leq\log|t|+1.
\]
In particular, if\, $t\in M_{a}$, then $|t|\leq4^{d}\max\{|a|,4\}^{d}$.
\end{prop}

\begin{proof}
For each $t\in\C$, let
\[
\lambda_{t}(z)=\lim_{n\to\infty}\frac{1}{d^{n}}\log^{+}|f_{t}^{n}(z)|.
\]
There exists an analytic homeomorphism $\phi_{t}$, defining the \textit{B\"ottcher coordinate} near $\infty$, which satisfies $\phi_{t}(f_{t}(z))=\phi_{t}(z)^{d}$ and $\lambda_{t}(z)=\log|\phi_{t}(z)|$. The map $\phi_{t}$ sends the domain
\[
V_{t}=\{z\in\C:\lambda_{t}(z)>\lambda_{t}(0)\}
\]
biholomorphically to $U_{R_{t}}$ with $R_{t}=e^{\lambda_{t}(0)}$.
By Theorem~\ref{thm3.1}, we have $V_{t}=\phi_{t}^{-1}(U_{R_{t}})\supseteq U_{2R_{t}}$.
By the proof of \cite[Lemma 3.2]{BD1}, we have $|t|\leq2^{d}R_{t}^{d}$,
and if $t$ is large enough such that $R_{t}^{d}-4R_{t}>2|a|^{d}$,
then $a^{d}+t\in U_{2R_{t}}$. If $|t|>4^{d}\max\{|a|,4\}^{d}$, then
we have
\[
R_{t}\geq\frac{1}{2}|t|^{1/d}>2\max\{|a|,4\}
\]
and
\begin{align*}
R_{t}^{d}-4R_{t} & =R_{t}\left(R_{t}^{d-1}-4\right)>2\max\{|a|,4\}\left(2^{d-1}\max\{|a|,4\}^{d-1}-4\right)\\
 & =2\max\{|a|,4\}^{d}+\left(\left(2^{d}-2\right)\max\{|a|,4\}^{d-1}-8\right)\max\{|a|,4\}\geq2|a|^{d}.
\end{align*}
Applying Theorem~\ref{thm3.1} to $\phi_{t}$, we get
\[
\left|\phi_{t}\left(a^{d}+t\right)\right|\leq2\left|a^{d}+t\right|\leq2\left(|a|^{d}+|t|\right)\leq\frac{17}{8}|t|,
\]
where the last inequality follows from $|t|>4^{d}\max\{|a|,4\}^{d}\geq16|a|^{d}$.
Applying Theorem~\ref{thm3.1} to $\phi_{t}^{-1}$, we get
\[
\left|\phi_{t}\left(a^{d}+t\right)\right|\geq\frac{1}{2}\left|\phi_{t}^{-1}\left(\phi_{t}\left(a^{d}+t\right)\right)\right|=\frac{1}{2}\left|a^{d}+t\right|\geq\frac{1}{2}\left(|t|-|a|^{d}\right)\geq\frac{15}{32}|t|.
\]
By \cite[Proposition 3.3]{BD1}, we have
\[
\log|t|-1\leq g_{a}(t)=\log\left|\phi_{t}\left(a^{d}+t\right)\right|\leq\log|t|+1.
\]
In particular, by Theorem~\ref{thm2.1}, if $|t|>4^{d}\max\{|a|,4\}^{d}$,
then $g_{a}(t)>0$ and $t\notin M_{a}$.
\end{proof}

By definition, $O_{f_{t}}(a)$ is bounded if $t\in M_{a}$. The following
result shows that $O_{f_{t}}(a)$ is in fact uniformly bounded for
any $t\in M_{a}$.

\begin{prop}\label{prop3.3}
If\, $t\in M_{a}$, then $|f_{t}^{n}(a)|<8\max\{|a|,4\}$
for any $n\geq1$.
\end{prop}

\begin{proof}
By Proposition~\ref{prop3.2}, if $t\in M_{a}$, then $|t|\leq4^{d}\max\{|a|,4\}^{d}$.
Suppose there exists an $n\geq1$ such that $|f_{t}^{n}(a)|\geq8\max\{|a|,4\}$;
then
\[
\left|f_{t}^{n+1}(a)\right|\geq\left|f_{t}^{n}(a)\right|^{d}-|t|\geq\left(1-\frac{1}{2^{d}}\right)\left|f_{t}^{n}(a)\right|^{d}\geq24\left|f_{t}^{n}(a)\right|.
\]
By induction, $|f_{t}^{n+k}(a)|\geq24^{k}|f_{t}^{n}(a)|\to\infty$
as $k\to\infty$, and we have a contradiction.
\end{proof}

The explicit bound of $M_{a}$ given in Proposition~\ref{prop3.2}
can be improved as follows.

\begin{prop}\label{prop3.4}
  If\, $t\in M_{a}$, then $|t|<3\max\{|a|,4\}^{d}$.
\end{prop}

\begin{proof}
If $t\in M_{a}$, then by Proposition~\ref{prop3.3},
\begin{align*}
|t| & \leq\left|a^{d}+t\right|+|a|^{d}<8\max\{|a|,4\}+|a|^{d}\\
 & \leq\left(\frac{8}{\max\{|a|,4\}^{d-1}}+1\right)\max\{|a|,4\}^{d}\leq3\max\{|a|,4\}^{d}.\qedhere
\end{align*}
\end{proof}

As \cite[Section 3.3]{DKY2}, we estimate the upper bound for $g_{a}(s)$
when $s$ is close to $M_{a}$.

\begin{prop}\label{prop3.5}
  If\, $t\in M_{a}$ and $|s-t|\leq\max\{|a|,4\}$, then for any $n\geq1$,
\[
\left|f_{s}^{n}(a)-f_{t}^{n}(a)\right|\leq|s-t|(18\max\{|a|,4\})^{d^{n-1}-1}.
\]
\end{prop}

\begin{proof}
Let $A_{n}=(18\max\{|a|,4\})^{d^{n-1}-1}$ for any $n\geq1$. We prove
the assertion by induction. 
The statement is true for $n=1$ because $|f_{s}(a)-f_{t}(a)|=|s-t|$. 
Assume the statement is true for
some $n\geq1$; then
\allowdisplaybreaks
\begin{align*}
\left|f_{s}^{n+1}(a)-f_{t}^{n+1}(a)\right|
=\; & \left|f_{s}^{n}(a)^{d}+s-f_{t}^{n}(a)^{d}-t\right|\\
\leq\; & \left|f_{s}^{n}(a)^{d}-f_{t}^{n}(a)^{d}\right|+|s-t|\\
=\; & {\textstyle \prod_{i=0}^{d-1}}\left|f_{s}^{n}(a)-f_{t}^{n}(a)+\left(1-\zeta_{d}^{i}\right)f_{t}^{n}(a)\right|+|s-t|\\
\leq\; & |f_{s}^{n}(a)-f_{t}^{n}(a)|(|f_{s}^{n}(a)-f_{t}^{n}(a)|+2|f_{t}^{n}(a)|)^{d-1}+|s-t|\\
\leq\; & |s-t|A_{n}(A_{n}\max\{|a|,4\}+16\max\{|a|,4\})^{d-1}+|s-t|\\
 & \text{(by the assumption, the induction hypothesis, and Proposition~\ref{prop3.3})}\\
=\; & |s-t|\left(A_{n}(A_{n}+16)^{d-1}\max\{|a|,4\}^{d-1}+1\right)\\
\leq\; & |s-t|\left(A_{n}(A_{n}+16A_{n})^{d-1}\max\{|a|,4\}^{d-1}+\max\{|a|,4\}^{d-1}\right)\\
=\; & |s-t|\left(17^{d-1}A_{n}^{d}+1\right)\max\{|a|,4\}^{d-1}\\
\leq\; & |s-t|\,18^{d-1}A_{n}^{d}\max\{|a|,4\}^{d-1}\\
=\; & |s-t|\,A_{n+1}.
\end{align*}
This completes the inductive step and hence the proof.
\end{proof}

\begin{prop}\label{prop3.6}
  If\, $t\in M_{a}$ and
\[
|s-t|\leq\frac{1}{18}(18\max\{|a|,4\})^{2-d^{n-1}}
\]
for some $n\geq1$, then
\[
g_{a}(s)\leq\frac{1}{d^{n-1}}\log(10\max\{|a|,4\}).
\]
\end{prop}

\begin{proof}
Let $A_{n}=(18\max\{|a|,4\})^{d^{n-1}-1}$ for any $n\geq1$. By Propositions
\ref{prop3.3} and~\ref{prop3.5},
\[
|f_{s}^{n}(a)|\leq|f_{s}^{n}(a)-f_{t}^{n}(a)|+|f_{t}^{n}(a)|\leq|s-t|A_{n}+|f_{t}^{n}(a)|\leq9\max\{|a|,4\}.
\]
By Proposition~\ref{prop3.4},
\[
\left|f_{s}^{n+1}(a)\right|\leq|f_{s}^{n}(a)|^{d}+|t|+|s-t|\leq\left(9^{d}+4\right)\max\{|a|,4\}^{d}.
\]
Let $p(z)=z^{d}+4$; then by induction and Lemma~\ref{lem3.7},
\[
\left|f_{s}^{n+k}(a)\right|\leq p^{k}(9)\max\{|a|,4\}^{d^{k}}\leq(10\max\{|a|,4\})^{d^{k}}
\]
for any $k\geq1$. Therefore, by Theorem~\ref{thm2.1},
\[
g_{a}(s)=\lim_{k\to\infty}\frac{1}{d^{n+k-1}}\log^{+}\left|f_{s}^{n+k}(a)\right|\leq\frac{1}{d^{n-1}}\log(10\max\{|a|,4\}).\qedhere
\]
\end{proof}

\begin{lem}\label{lem3.7}
Let $p(z)=z^{d}+4$ for some $d\geq2$. Then $p^{n}(9)\leq10^{d^{n}}$
for any $n\geq1$.
\end{lem}

\begin{proof}
We prove $q(n)=10^{d^{n}}-p^{n}(9)\geq1$ for any $n\geq0$ by induction.
It is clear that $q(0)=1$. Assume $q(n)\geq1$ for some $n\geq0$;
then
\begin{align*}
q(n+1) & =\left(10^{d^{n}}\right)^{d}-p^{n}(9)^{d}-4=q(n){\textstyle \sum_{i=0}^{d-1}}\left(10^{d^{n}}\right)^{i}p^{n}(9)^{d-1-i}-4\\
 & \geq q(n)10^{d^{n}(d-1)}-4\geq1.
\end{align*}
This completes the inductive step and hence the proof.
\end{proof}

\subsection{Non-Archimedean estimates}\label{sec3.2}

In this section, we assume that $K$ is a number field, $a\in K$,
and $v\in M_{K}^{0}$. Because $v$ is fixed, we write $|\cdot|$,
$\delta$, $M_{a}$, and $g_{a}$ for $|\cdot|_{v}$, $\delta_{v}$,
$M_{a,v}$, and $g_{a,v}$.

Propositions~\ref{prop3.8} and~\ref{prop3.9} can be seen as the
non-Archimedean versions of Propositions~\ref{prop3.2},~\ref{prop3.3},
and~\ref{prop3.4}.

\begin{prop}\label{prop3.8}
  Assume that $|a|\leq1$. For any $t\in\C_{v}$,
we have
\[
g_{a}(t)=\log^{+}|t|.
\]
In particular, $M_{a}=\mathcal{D}(0,1)$, the closed Berkovich unit
disk.
\end{prop}

\begin{proof}
Let $t\in\C_{v}$; then for any $n\geq1$,
\[
|f_{t}^{n}(a)|\begin{cases}
\leq1 & \text{if }|t|\leq1,\\
=|t|^{d^{n-1}} & \text{if }|t|>1.
\end{cases}
\]
By Theorem~\ref{thm2.1},
\[
g_{a}(t)=\lim_{n\to\infty}\frac{1}{d^{n-1}}\log^{+}|f_{t}^{n}(a)|=\log^{+}|t|.
\]
Therefore, $M_{a}\cap\C_{v}=\overline{D}(0,1)$ and $M_{a}=\overline{M_{a}\cap\C_{v}}=\mathcal{D}(0,1)$.
\end{proof}

\begin{prop}\label{prop3.9}
Assume that $|a|>1$. For any $t\in\C_{v}$, we have
\[
g_{a}(t)\begin{cases}
=d\log|a| & \text{if }\,|t|<|a|^{d},\\
\leq d\log|a| & \text{if }\,|t|=|a|^{d},\\
=\log|t| & \text{if }\,|t|>|a|^{d}.
\end{cases}
\]
If\, $t\in M_{a}\cap\C_{v}$, then $|t|=|a|^{d}$ and $|f_{t}^{n}(a)|=|a|$
for any $n\geq1$.
\end{prop}

\begin{proof}
Let $t\in\C_{v}$; then for any $n\geq1$,
\[
|f_{t}^{n}(a)|\begin{cases}
=|a|^{d^{n}} & \text{if }|t|<|a|^{d},\\
\leq|a|^{d^{n}} & \text{if }|t|=|a|^{d},\\
=|t|^{d^{n-1}} & \text{if }|t|>|a|^{d}.
\end{cases}
\]
By Theorem~\ref{thm2.1},
\[
g_{a}(t)=\lim_{n\to\infty}\frac{1}{d^{n-1}}\log^{+}|f_{t}^{n}(a)|\begin{cases}
=d\log|a| & \text{if }|t|<|a|^{d},\\
\leq d\log|a| & \text{if }|t|=|a|^{d},\\
=\log|t| & \text{if }|t|>|a|^{d}.
\end{cases}
\]
Therefore, if $t\in M_{a}\cap\C_{v}$, then $|t|=|a|^{d}$
and $t\in M_{f_{t}^{n}(a)}\cap\C_{v}$ for any $n\geq1$.
By the same reasoning, we have $|t|=|f_{t}^{n}(a)|^{d}$ and $|f_{t}^{n}(a)|=|a|$.
\end{proof}

The following result will not be used until Section~\ref{sec4}. We
include it here because it is similar to Proposition~\ref{prop3.9}.

\begin{prop}\label{prop3.10}
If\, $|a|>1$, $t\in\C_{v}$, and $|f_{t}^{n}(a)|<|a|^{d}$
for some $n\geq1$, then $|t|=|a|^{d}$ and $|f_{t}^{k}(a)|=|a|$
for any $1\leq k\leq n-1$.
\end{prop}

\begin{proof}
For the first assertion, suppose $|t|\neq|a|^{d}$; then
\[
|f_{t}^{n}(a)|=\max\left\{|t|,|a|^{d}\right\}^{d^{n-1}}\geq|a|^{d}.
\]
For the second assertion, suppose $|f_{t}^{k}(a)|\neq|a|$ for some
$1\leq k\leq n-1$; then
\[
|f_{t}^{n}(a)|=\max\left\{|t|,\left|f_{t}^{k}(a)\right|^{d}\right\}^{d^{n-k-1}}\geq|t|=|a|^{d}.\qedhere
\]
\end{proof}

Propositions~\ref{prop3.11} and~\ref{prop3.12} can be seen as the
non-Archimedean versions of Propositions~\ref{prop3.5} and~\ref{prop3.6}.

\begin{prop}\label{prop3.11}
If\, $|a|>1$, $t\in M_{a}\cap\C_{v}$, $s\in\C_{v}$, and $|s-t|\leq|a|$, then 
for any $n\geq1$, 
\[
|f_{s}^{n}(a)-f_{t}^{n}(a)|\leq|s-t|\cdot|a|^{d^{n-1}-1}. 
\]
\end{prop}

\begin{proof}
We prove the assertion by induction. Because $|f_{s}(a)-f_{t}(a)|=|s-t|$,
the statement is true for $n=1$. Assume the statement is true for
some $n\geq1$; then
\allowdisplaybreaks
\begin{align*}
 \left|f_{s}^{n+1}(a)-f_{t}^{n+1}(a)\right|
=\; & \left|f_{s}^{n}(a)^{d}+s-f_{t}^{n}(a)^{d}-t\right|\\
\leq\; & \max\left\{\left|f_{s}^{n}(a)^{d}-f_{t}^{n}(a)^{d}\right|,|s-t|\right\}\\
=\; & \max\left\{{\textstyle \prod_{i=0}^{d-1}}\left|f_{s}^{n}(a)-f_{t}^{n}(a)+(1-\zeta_{d}^{i})f_{t}^{n}(a)\right|,|s-t|\right\}\\
\leq\; & \max\left\{|f_{s}^{n}(a)-f_{t}^{n}(a)|\max\{|f_{s}^{n}(a)-f_{t}^{n}(a)|,|f_{t}^{n}(a)|\}^{d-1},|s-t|\right\}\\
\leq\; & \max\left\{|s-t|\cdot|a|^{d^{n-1}-1}\max\left\{|a|^{d^{n-1}},|a|\right\}^{d-1},|s-t|\right\}\\
 & \text{(by the assumption, the induction hypothesis, and Proposition~\ref{prop3.9})}\\
=\; & |s-t|\cdot|a|^{d^{n}-1}.
\end{align*}
This completes the inductive step and hence the proof.
\end{proof}

\begin{prop}
\label{prop3.12}If $t\in M_{a}\cap\C_{v}$, $s\in\A_{v}^{1,\an}$,
and
\[
\delta(s,t)\leq\max\{|a|,1\}^{2-d^{n-1}}
\]
for some $n\geq1$, then
\[
g_{a}(s)\leq\frac{1}{d^{n-1}}\log^{+}|a|.
\]
\end{prop}

\begin{proof}
Let $A_{n}=\max\{|a|,1\}^{2-d^{n-1}}$ for any $n\geq1$. By Theorem
\ref{thm2.1}, the function $g_{a}$ is continuous. Since $\overline{D}(t,A_{n})$
is dense in $\mathcal{D}(t,A_{n})$, it suffices to prove the assertion
for $s\in\overline{D}(t,A_{n})$. If $|a|\leq1$, then by Proposition~\ref{prop3.8}, we have $s\in M_{a}$ and $g_{a}(s)=0$. If $|a|>1$, then
by Propositions~\ref{prop3.9} and~\ref{prop3.11},
\[
|f_{s}^{n}(a)|\leq\max\{|f_{s}^{n}(a)-f_{t}^{n}(a)|,|f_{t}^{n}(a)|\}\leq|a|
\]
and
\[
\left|f_{s}^{n+1}(a)\right|\leq\max\left\{|f_{s}^{n}(a)|^{d},|t|,|s-t|\right\}\leq|a|^{d}.
\]
By induction, $|f_{s}^{n+k}(a)|\leq|a|^{d^{k}}$ for any $k\geq1$.
Therefore, by Theorem~\ref{thm2.1},
\[
g_{a}(s)=\lim_{k\to\infty}\frac{1}{d^{n+k-1}}\log^{+}\left|f_{s}^{n+k}(a)\right|\leq\frac{1}{d^{n-1}}\log|a|.\qedhere
\]
\end{proof}

\subsection{Proof of Theorem \texorpdfstring{\ref{thm1.2}}{1.2}}\label{sec3.3}

Now we are ready to give the proof of Theorem~\ref{thm1.2}. Given
$a,b\in\bar{\Q}$, let $K$ be a number field such that $a,b\in K$.
As in \cite[Section 9]{DKY2}, we will apply \cite[Theorem 1]{Fi} in
the following way:
\[
\left\langle \boldsymbol{\mu}_{a},\boldsymbol{\mu}_{b}\right\rangle ^{1/2}\leq\left\langle \boldsymbol{\mu}_{a},[S]_{\boldsymbol{\tau}}\right\rangle ^{1/2}+\left\langle \boldsymbol{\mu}_{b},[S]_{\boldsymbol{\tau}}\right\rangle ^{1/2},
\]
where $[S]_{\boldsymbol{\tau}}$ is an adelic measure to be described
below.

Let $S$ be a finite, non-empty, $\Gal(\bar{K}/K)$-invariant
subset of $\bar{K}$, and let $[S]$ be the probability measure supported
equally on the elements of $S$. We call $\boldsymbol{\tau}=\{\tau_{v}\}_{v\in M_{K}}$
an \textit{adelic radius} if $\tau_{v}>0$ for any $v\in M_{K}$
and $\tau_{v}=1$ for all but finitely many $v\in M_{K}$. Combining
the ideas of \cite[Section 4.6]{FRL} and \cite[Section 2.1]{FP},
we define the \textit{regularization} $[S]_{\boldsymbol{\tau}}=\{[S]_{\tau_{v}}\}_{v\in M_{K}}$
by
\[
[S]_{\tau_{v}}=\frac{1}{|S|}\sum_{s\in S}m_{s,\tau_{v},v},
\]
where $m_{s,\tau_{v},v}$ is defined in Section~\ref{sec2.4}.

The following result and the proof of Theorem~\ref{thm1.2} are adapted
from \cite[Section 9]{DKY2}.

\begin{thm}\label{thm3.13}
Let $K$ be a number field such that $a,b\in K$.
Then
\[
\left\langle \boldsymbol{\mu}_{a},\boldsymbol{\mu}_{b}\right\rangle ^{1/2}\leq\sum_{i=a,b}\left(\sum_{v\in M_{K}}n_{v}\left(-(\mu_{i,v},[S]_{\tau_{v}})_{v}-\frac{\log\tau_{v}}{2|S|}\right)\right)^{1/2}
\]
for any finite, non-empty, $\Gal(\bar{K}/K)$-invariant subset
$S$ of $\bar{K}$ and any adelic radius $\boldsymbol{\tau}$.
\end{thm}

\begin{proof}
By Theorem~\ref{thm2.1}, we have $(\mu_{i,v},\mu_{i,v})_{v}=-\log\gamma_{v}(M_{i,v})=0$
for $i=a,b$ and any $v\in M_{K}$. Then the proof is identical to
that of \cite[Lemma 9.2]{DKY2}.
\end{proof}

\begin{proof}[Proof of Theorem~\ref{thm1.2}]
Let $K$ be a number field such
that $a,b\in K$, and let $n\geq-1$ be the integer such that $d^{-n-1}\leq\varepsilon'=\varepsilon/4<d^{-n}$.
For each $v\in M_{K}^{\infty}$, take
\begin{align*}
\tau_{v} & =(18\max\{|a|_{v},|b|_{v},4\})^{1-d/\varepsilon'}\leq(18\max\{|a|_{v},4\})^{1-d/\varepsilon'}\\
 & \leq\frac{1}{18}(18\max\{|a|_{v},4\})^{2-d/\varepsilon'}\leq\frac{1}{18}(18\max\{|a|_{v},4\})^{2-d^{n+1}}.
\end{align*}
By Proposition~\ref{prop3.6}, if $|s-t|_{v}=\tau_{v}$ for some $t\in S_{a,b}$,
then
\[
g_{a,v}(s)\leq\frac{1}{d^{n+1}}\log(10\max\{|a|_{v},4\})\leq\varepsilon'\log(18\max\{|a|_{v},|b|_{v},4\}).
\]
Therefore,
\begin{align*}
 -\left(\mu_{a,v},[S_{a,b}]_{\tau_{v}}\right)_{v}-\frac{\log\tau_{v}}{2|S_{a,b}|}=\;&\int g_{a,v}d[S_{a,b}]_{\tau_{v}}-\frac{\log\tau_{v}}{2|S_{a,b}|}\\
\leq\; & \varepsilon'\log(18\max\{|a|_{v},|b|_{v},4\})-\frac{(1-d/\varepsilon')\log(18\max\{|a|_{v},|b|_{v},4\})}{2|S_{a,b}|}\\
\leq\; & \left(\varepsilon'+\frac{d/\varepsilon'-1}{2|S_{a,b}|}\right)\log(18\max\{4|a|_{v},4|b|_{v},4\})\\
\leq\; & \left(\varepsilon'+\frac{d/\varepsilon'-1}{2|S_{a,b}|}\right)(\log\max\{|a|_{v},|b|_{v},1\}+5).
\end{align*}
For each $v\in M_{K}^{0}$, take
\begin{align*}
\tau_{v} & =\max\{|a|_{v},|b|_{v},1\}^{1-d/\varepsilon'}\leq\max\{|a|_{v},1\}^{1-d/\varepsilon'}\\
 & \leq\max\{|a|_{v},1\}^{2-d/\varepsilon'}\leq\max\{|a|_{v},1\}^{2-d^{n+1}}.
\end{align*}
By Proposition~\ref{prop3.12}, if $\delta_{v}(s,t)=\tau_{v}$ for
some $t\in S_{a,b}$, then
\[
g_{a,v}(s)\leq\frac{1}{d^{n+1}}\log^{+}|a|_{v}\leq\varepsilon'\log\max\{|a|_{v},|b|_{v},1\}.
\]
Therefore,
\begin{align*}
  -\left(\mu_{a,v},[S_{a,b}]_{\tau_{v}}\right)_{v}-\frac{\log\tau_{v}}{2|S_{a,b}|}=\;&\int g_{a,v}d[S_{a,b}]_{\tau_{v}}-\frac{\log\tau_{v}}{2|S_{a,b}|}\\
\leq\; & \varepsilon'\log\max\{|a|_{v},|b|_{v},1\}-\frac{(1-d/\varepsilon')\log\max\{|a|_{v},|b|_{v},1\}}{2|S_{a,b}|}\\
=\; & \left(\varepsilon'+\frac{d/\varepsilon'-1}{2|S_{a,b}|}\right)\log\max\{|a|_{v},|b|_{v},1\}.
\end{align*}
Summing over all $v\in M_{K}$, we get
\[
\left\langle \boldsymbol{\mu}_{a},\boldsymbol{\mu}_{b}\right\rangle \leq4\left(\varepsilon'+\frac{d/\varepsilon'-1}{2|S_{a,b}|}\right)(h(a,b)+5)=\left(\varepsilon+\frac{8d/\varepsilon-2}{|S_{a,b}|}\right)(h(a,b)+5).\qedhere
\]
\end{proof}

\section{Lower bounds}\label{sec4}

The purpose of this section is to give the proofs of Theorems~\ref{thm1.3}
and~\ref{thm1.4}. In order to do so, we estimate the lower bounds
for
\[
\left\langle \mu_{a,v},\mu_{b,v}\right\rangle _{v}=\int_{\A_{v}^{1,\an}}g_{a,v}d\mu_{b,v}.
\]
In Section~\ref{sec4.1}, we review two equidistribution theorems
for later usage. As always, we work with $v\in M_{K}^{\infty}$ and
$v\in M_{K}^{0}$ separately. For some technical reasons, we also
work with $d=2$ and $d>2$ separately when $v\in M_{K}^{\infty}$.
The local estimates obtained in Sections~\ref{sec4.2},~\ref{sec4.3},
and~\ref{sec4.4} are gathered in Section~\ref{sec4.5} to complete
the proofs of Theorems~\ref{thm1.3} and~\ref{thm1.4}.

\subsection{Equidistribution theorems}\label{sec4.1}

In this section, we assume that $K$ is a number field and $a\in K$.
We give two equidistribution theorems. Theorem~\ref{thm4.1} will
be used in Section~\ref{sec4.4}, and Theorem~\ref{thm4.3} will be
used in Sections~\ref{sec4.2} and~\ref{sec4.3}.

We call $\boldsymbol{E}=\{E_{v}\}_{v\in M_{K}}$ an \textit{adelic
compact set} if $E_{v}$ is a non-empty compact subset of $\A_{v}^{1,\an}$
for any $v\in M_{K}$ and $E_{v}=\mathcal{D}(0,1)$ for all but finitely
many $v\in M_{K}$.

\begin{thm}[\textit{cf.} {\cite[Theorem 7.52]{BR2}}]\label{thm4.1}
Let $K$ be a number field, and let $\boldsymbol{E}$ be an adelic compact set with
\[
\gamma(\boldsymbol{E})=\prod_{v\in M_{K}}\gamma_{v}(E_{v})^{n_{v}}=1.
\]
Suppose $S_{n}$ is a sequence of finite, non-empty, $\Gal(\bar{K}/K)$-invariant
subsets of\, $\bar{K}$ such that $|S_{n}|\to\infty$ and
\[
h_{\boldsymbol{E}}(S_{n})=\sum_{v\in M_{K}}n_{v}\left(\frac{1}{|S_{n}|}\sum_{z\in S_{n}}g_{E_{v}}(z)\right)\longrightarrow 0.
\]
Fix a $v\in M_{K}$ and, for any $n\geq1$, let $[S_{n}]$ be the probability
measure on $\A_{v}^{1,\an}$ supported equally on
the elements of\,~$S_{n}$. Then the sequence of measures $[S_{n}]$
converges weakly to $\mu_{E_{v}}$ on $\A_{v}^{1,\an}$.
\end{thm}

For  related results, see also \cite{Bi,BR1,CL1,FRL} and \cite[Theorem 10.24]{BR2}.

From Theorem~\ref{thm2.1} and Proposition~\ref{prop3.8}, we know
that $\boldsymbol{M}_{a}=\{M_{a,v}\}_{v\in M_{K}}$ is an adelic compact
set with $\gamma(\boldsymbol{M}_{a})=1$. Let $S_{n}$ be the set
of all roots of $f_{T}^{n}(a)$. The following result shows that $\mu_{a,v}$
can be approximated by $[S_{n}]$ if we assume that all roots of $f_{T}^{n}(a)$
are simple.

\begin{prop}\label{prop4.2}
Assume that, for any $n\geq1$, all roots of $f_{T}^{n}(a)$
are simple. Let $S_{n}$ be the set of all roots of $f_{T}^{n}(a)$.
Then for any $v\in M_{K}$, $[S_{n}]$ converges weakly to $\mu_{a,v}$
on $\A_{v}^{1,\an}$.
\end{prop}

\begin{proof}
By the assumption, we have $|S_{n}|=d^{n-1}\to\infty$ as $n\to\infty$.
To apply Theorem~\ref{thm4.1}, it remains to show that $h_{\boldsymbol{M}_{a}}(S_{n})\to0$
as $n\to\infty$. If $t\in S_{n}$, then by Theorem~\ref{thm2.1},
\[
g_{a,v}(t)=\lim_{k\to\infty}\frac{1}{d^{n+k}}\log^{+}\left|f_{t}^{n+k+1}(a)\right|_{v}=\frac{1}{d^{n}}\lim_{k\to\infty}\frac{1}{d^{k}}\log^{+}\left|f_{t}^{k+1}(0)\right|_{v}=\frac{1}{d^{n}}g_{0,v}(t).
\]
Since $g_{0,v}$ is continuous, it suffices to show that for any $v\in M_{K}$
and any $t\in S_{n}$, $|t|_{v}$ is bounded by some constant $C_{a,v}$,
where $C_{a,v}$ is independent of $n$.
\begin{enumerate}[wide]
\item If $v\in M_{K}^{0}$ and $|a|_{v}\leq1$, then by Proposition
\ref{prop3.8}, we have $|t|_{v}\leq1$ for any $t\in S_{n}$.
\item If $v\in M_{K}^{0}$ and $|a|_{v}>1$, then by Proposition~\ref{prop3.10},
we have $|t|_{v}=|a|_{v}^{d}$ for any $t\in S_{n}$.
\item If $v\in M_{K}^{\infty}$, then by Proposition~\ref{prop3.4},
the roots of $f_{T}^{n}(a)=a$ are inside
\[
D_{a,v}=D(0,3\max\{|a|_{v},4\}^{d}).
\]
Once we can show that $|f_{t}^{n}(a)|_{v}>|a|_{v}$ for any $t\in\partial D_{a,v}$,
we can apply Rouch\'e's theorem to conclude $S_{n}\subseteq D_{a,v}$.
If $t\in\partial D_{a,v}$, then
\[
|f_{t}(a)|_{v}\geq|t|_{v}-|a|_{v}^{d}=3\max\{|a|_{v},4\}^{d}-|a|_{v}^{d}\geq2\max\{|a|_{v},4\}^{d}\geq8\max\{|a|_{v},4\}.
\]
Assume that $|f_{t}^{n}(a)|_{v}\geq8\max\{|a|_{v},4\}$ for some $n\geq1$; 
then
\[
\left|f_{t}^{n+1}(a)\right|_{v}\geq|f_{t}^{n}(a)|_{v}^{d}-|t|_{v}\geq\left(8^{d}-3\right)\max\{|a|_{v},4\}^{d}\geq8\max\{|a|_{v},4\}.
\]
By induction, $|f_{t}^{n}(a)|_{v}\geq8\max\{|a|_{v},4\}>|a|_{v}$
for any $n\geq1$.\hfill\qedhere
\end{enumerate}
\end{proof}

Proposition~\ref{prop4.10} shows that if $|a|_{v}>|d|_{v}^{-2/(d-1)}$
for some $v\in M_{K}^{0}$, then all roots of $f_{T}^{n}(a)$ are
simple. Since we lack a similar result for $v\in M_{K}^{\infty}$,
we need an equidistribution theorem taking the multiplicity into account.

\begin{thm}\label{thm4.3}
Fix a $v\in M_{K}^{\infty}$ and a sequence $0\leq k(n)<n$.
Let $\delta_{n}$ be the discrete probability measure on $\C$
weighted by multiplicity on the roots of $f_{T}^{n}(a)=f_{T}^{k(n)}(a)$.
Then the sequence of measures $\delta_{n}$ converges weakly to $\mu_{a,v}$
on $\C$.
\end{thm}

\begin{proof}
See \cite[Section 4.3]{BD1} and the proof of \cite[Theorem 1]{DF}.
\end{proof}

\subsection{Archimedean estimates for \texorpdfstring{$\boldsymbol{d=2}$}{$d=2$}}\label{sec4.2}

In this section, we assume that $K$ is a number field, $a,b\in K$,
$v\in M_{K}^{\infty}$, and $d=2$. Because $v$ is fixed, we write
$|\cdot|$, $M_{a}$, $\mu_{a}$, and $g_{a}$ for $|\cdot|_{v}$,
$M_{a,v}$, $\mu_{a,v}$, and $g_{a,v}$.

As \cite[Section 3.2]{DKY2}, we first cover $M_{a}$ by some disjoint
open disks around the roots of $f_{T}^{n}(a)=f_{T}(a)$ with $n=2,3$,
and then we estimate the lower bound for $g_{a}(t)$ when $t$ is outside
the cover. Let
\[
\alpha_{1}(a)=-a^{2}-a,\quad\alpha_{2}(a)=-a^{2}+a,\quad\alpha_{3}(a)=-a^{2}-a-1,\quad\alpha_{4}(a)=-a^{2}+a-1
\]
be the roots of $f_{T}^{3}(a)=f_{T}(a)$. Note that $\alpha_{1}(a)$
and $\alpha_{2}(a)$ are also roots of $f_{T}^{2}(a)=f_{T}(a)$.

\begin{prop}\label{prop4.4}
  Assume that $|a|\geq28$. Then we have
\begin{enumerate}
\item\label{p4.4-1} $M_{a}\subseteq\cup_{i=1}^{2}D(\alpha_{i}(a),10)$; 
\item\label{p4.4-2} if $t\notin\cup_{i=1}^{2}D(\alpha_{i}(a),10)$, then
\[
g_{a}(t)\geq\frac{1}{2}\log(13|a|);
\]
\item\label{p4.4-3} $\mu_{a}(D(\alpha_{i}(a),10))=1/2$ for $i=1,2$.
\end{enumerate}
\end{prop}

\begin{proof}
  For simplicity, let $\alpha_{i}=\alpha_{i}(a)$ and $D_{i}=D(\alpha_{i},10)$.

\eqref{p4.4-1}~ If $t\in\partial D_{1}$, then
\begin{align*}
|t-\alpha_{1}| & =10,\\
|t-\alpha_{2}| & \geq|\alpha_{1}-\alpha_{2}|-|t-\alpha_{1}|=2|a|-10.
\end{align*}
When $|a|\geq28$, we have
\[
|f_{t}^{2}(a)-f_{t}(a)|\geq20|a|-100\geq16|a|.
\]
The same reasoning also works for $t\in\partial D_{2}$. By Rouch\'e's
theorem, for any $c$ with $|c|<16|a|$, the equation $f_{T}^{2}(a)-f_{T}(a)=c$
has exactly one root in each $D_{i}$. By Proposition~\ref{prop3.3},
\[
\left|f_{s}^{2}(a)-f_{s}(a)\right|\leq\left|f_{s}^{2}(a)\right|+|f_{s}(a)|<16|a|
\]
for any $s\in M_{a}$, so we have $M_{a}\subseteq\cup_{i=1}^{2}D_{i}$.

\eqref{p4.4-2}~ If $t\in\partial D_{i}$ for some $i=1,2$, then
\begin{align*}
|t| & \leq|t-\alpha_{i}|+|\alpha_{i}|\leq2|a|^{2},\\
|f_{t}(a)| & \leq|t-\alpha_{i}|+|\alpha_{i}+a^{2}|\leq2|a|
\end{align*}
and
\[
\left|f_{t}^{2}(a)\right|\geq\left|f_{t}^{2}(a)-f_{t}(a)\right|-|f_{t}(a)|\geq14|a|.
\]
Let $p(z)=z^{2}-2$; then by induction and an argument similar to that of Lemma
\ref{lem3.7},
\[
\left|f_{t}^{n+2}(a)\right|\geq p^{n}(14)|a|^{2^{n}}\geq(13|a|)^{2^{n}}
\]
for any $n\geq1$. Therefore, by Theorem~\ref{thm2.1},
\[
g_{a}(t)=\lim_{n\to\infty}\frac{1}{2^{n+1}}\log^{+}\left|f_{t}^{n+2}(a)\right|\geq\frac{1}{2}\log(13|a|).
\]
Since $g_{a}$ is harmonic on $\C\backslash M_{a}$, this
is true for any $t\notin\cup_{i=1}^{2}D_{i}$.

\eqref{p4.4-3}~  If $t\in\partial D_{i}$ for some $i=1,2$, then
\[
|f_{t}(a)\pm a|\leq|f_{t}(a)|+|a|\leq3|a|
\]
and
\[
\left|f_{t}^{n}(a)-f_{t}(a)\right|\geq\left|f_{t}^{n}(a)\right|-|f_{t}(a)|\geq(13|a|)^{2^{n-2}}-2|a|>3|a|
\]
for any $n\geq2$. By Rouch\'e's theorem, $f_{T}^{n}(a)-f_{T}(a)$
and $f_{T}^{n}(a)\pm a$ have the same number of roots in $D_{i}$.
Since
\[
f_{T}^{n+1}(a)-f_{T}(a)=\left(f_{T}^{n}(a)+a\right)\left(f_{T}^{n}(a)-a\right),
\]
by induction each of $f_{T}^{n}(a)-f_{T}(a)$ and $f_{T}^{n}(a)\pm a$
has $2^{n-2}$ roots in $D_{i}$. Then the conclusion follows from
Theorem~\ref{thm4.3}.
\end{proof}

\begin{prop}\label{prop4.5}
Assume that $|a|\geq28$. Then we have
\begin{enumerate}
\item $M_{a}\subseteq\cup_{i=1}^{4}D(\alpha_{i}(a),5/|a|)$;
\item if\, $t\notin\cup_{i=1}^{4}D(\alpha_{i}(a),5/|a|)$, then
\[
g_{a}(t)\geq\frac{1}{4}\log(13|a|); 
\]
\item $\mu_{a}(D(\alpha_{i}(a),5/|a|))=1/4$ for $1\leq i\leq4$.
\end{enumerate}
\end{prop}

\begin{proof}
For simplicity, let $\alpha_{i}=\alpha_{i}(a)$ and $D_{i}=D(\alpha_{i},5/|a|)$.
If $t\in\partial D_{1}$, then
\begin{align*}
|t-\alpha_{1}| & =5/|a|,\\
|t-\alpha_{2}| & \geq|\alpha_{1}-\alpha_{2}|-|t-\alpha_{1}|=2|a|-5/|a|,\\
|t-\alpha_{3}| & \geq|\alpha_{1}-\alpha_{3}|-|t-\alpha_{1}|=1-5/|a|,\\
|t-\alpha_{4}| & \geq|\alpha_{1}-\alpha_{4}|-|t-\alpha_{1}|\geq2|a|-1-5/|a|.
\end{align*}
When $|a|\geq28$, we have
\[
\left|f_{t}^{3}(a)-f_{t}(a)\right|\geq20|a|-110-\frac{50}{|a|}+\frac{525}{|a|^{2}}-\frac{625}{|a|^{4}}\geq16|a|.
\]
The rest of the proof is similar to that of Proposition~\ref{prop4.4}.
\end{proof}

As \cite[Theorem 4.1]{DKY2}, we estimate the complex Arakelov--Zhang
pairing as follows.

\begin{prop}\label{prop4.6}
For any $a,b\in K$, we have
\[
\int g_{a}d\mu_{b}\geq\frac{1}{8}\log^{+}\left|a^{2}-b^{2}\right|-\frac{1}{8}\log5000.
\]
Moreover, if\, $\max\{|a|,|b|\}\geq50$ and
\begin{equation}
\left|a^{2}-b^{2}\right|\geq\frac{11}{\max\{|a|,|b|\}},\label{eq4.1}
\end{equation}
then
\[
\int g_{a}d\mu_{b}\geq\frac{1}{16}\log\max\{|a|,|b|\}.
\]
\end{prop}

\begin{proof}
Without loss of generality, we assume that $|a|\geq|b|$. We prove
the second assertion in parts~\eqref{part1},~\eqref{part2},~\eqref{part3},  and we prove the first assertion in part~\eqref{part4}.

\begin{enumerate}[wide]
  \item\label{part1}
 \fbox{$|a|\geq50$ and $|b|\leq28$} If $t\in M_{b}$, then
by Proposition~\ref{prop3.4},
\[
|t|\leq3\max\{|b|,4\}^{2}\leq3\cdot28^{2}\leq50(50-2)\leq|a|(|a|-2)=|a|^{2}-2|a|
\]
and
\[
|f_{t}(a)|\geq|a|^{2}-|t|\geq2|a|.
\]
Let $p(z)=z^{2}-1$; then by induction and an argument similar to that of Lemma
\ref{lem3.7},
\[
\left|f_{t}^{n+1}(a)\right|\geq p^{n}(2)|a|^{2^{n}}\geq|a|^{2^{n}}
\]
for any $n\geq1$. Therefore, by Theorem~\ref{thm2.1},
\[
g_{a}(t)=\lim_{n\to\infty}\frac{1}{2^{n}}\log^{+}\left|f_{t}^{n+1}(a)\right|\geq\log|a|
\]
and
\[
\int g_{a}d\mu_{b}\geq\log|a|.
\]

\item\label{part2}
\fbox{$|a|\geq50$, $|b|\geq28$, and $|a^{2}-b^{2}|\geq20$}

\begin{claim*}
  Some $D(\alpha_{i}(b),10)$ is disjoint from $\cup_{i=1}^{2}D(\alpha_{i}(a),10)$.
\end{claim*}

Suppose the claim is false; then for each $i$, there exists a $k_{i}$ such that
\[
\left|\alpha_{k_{i}}(a)-\alpha_{i}(b)\right|<20.
\]
If $k_{1}=k_{2}$, then
\[
2|b|=|\alpha_{1}(b)-\alpha_{2}(b)|\leq\sum_{i=1}^{2}|\alpha_{k_{i}}(a)-\alpha_{i}(b)|<40,
\]
and we have a contradiction. If $k_{1}\neq k_{2}$, then
\[
2\left|a^{2}-b^{2}\right|=\left|\sum_{i=1}^{2}\alpha_{i}(a)-\sum_{i=1}^{2}\alpha_{i}(b)\right|\leq\sum_{i=1}^{2}\left|\alpha_{k_{i}}(a)-\alpha_{i}(b)\right|<40,
\]
and we also have a contradiction. Therefore, the claim is proved and, by Proposition
\ref{prop4.4},
\[
\int g_{a}d\mu_{b}\geq\frac{1}{2}\cdot\frac{1}{2}\log(13|a|)=\frac{1}{4}\log(13|a|).
\]

\item\label{part3}
  \fbox{$|a|\geq50$, $|b|\geq28$, and $11/|a|\leq|a^{2}-b^{2}|\leq20$}
  
\begin{claim*} Some $D(\alpha_{i}(b),5/|b|)$ is disjoint from $\cup_{i=1}^{4}D(\alpha_{i}(a),5/|a|)$.
\end{claim*}

Suppose the claim is false; then for each $i$, there exists a $k_{i}$ such that
\[
\left|\alpha_{k_{i}}(a)-\alpha_{i}(b)\right|<5/|a|+5/|b|.
\]
If $k_{i}=k_{j}$ for some $i\neq j$, then
\[
1\leq|\alpha_{i}(b)-\alpha_{j}(b)|\leq|\alpha_{k_{i}}(a)-\alpha_{i}(b)|+|\alpha_{k_{j}}(a)-\alpha_{j}(b)|<10/|a|+10/|b|,
\]
and we have a contradiction. If $k_{i}\neq k_{j}$ for any $i\neq j$, then
\[
4\left|a^{2}-b^{2}\right|=\left|\sum_{i=1}^{4}\alpha_{i}(a)-\sum_{i=1}^{4}\alpha_{i}(b)\right|\leq\sum_{i=1}^{4}|\alpha_{k_{i}}(a)-\alpha_{i}(b)|<20/|a|+20/|b|.
\]
This also gives a contradiction because
\[
\frac{|a|}{|b|}\leq\left(\left|\frac{a^{2}}{b^{2}}-1\right|+1\right)^{1/2}=\left(\frac{\left|a^{2}-b^{2}\right|}{|b|^{2}}+1\right)^{1/2}\leq\left(\frac{20}{28^{2}}+1\right)^{1/2}\leq\frac{6}{5}
\]
and
\[
\left|a^{2}-b^{2}\right|\geq11/|a|\geq5/|a|+5/|b|.
\]
Therefore, the claim is proved and, by Proposition~\ref{prop4.5},
\[
\int g_{a}d\mu_{b}\geq\frac{1}{4}\cdot\frac{1}{4}\log(13|a|)=\frac{1}{16}\log(13|a|).
\]

\item\label{part4}
Now we prove the first assertion. If $|a|\leq50$ or $|a^{2}-b^{2}|\leq20$,
then
\[
\log^{+}\left|a^{2}-b^{2}\right|\leq\log\max\left\{2|a|^{2},20\right\}\leq\log5000
\]
and
\[
\int g_{a}d\mu_{b}\geq0\geq\frac{1}{8}\log^{+}\left|a^{2}-b^{2}\right|-\frac{1}{8}\log5000.
\]
If $|a|\geq50$ and $|a^{2}-b^{2}|\geq20$, then by parts~\eqref{part1} and~\eqref{part2},
\[
\int g_{a}d\mu_{b}\geq\frac{1}{4}\log(13|a|)\geq\frac{1}{8}\log\left(2|a|^{2}\right)\geq\frac{1}{8}\log^{+}\left|a^{2}-b^{2}\right|.\qedhere
\]
\end{enumerate}
\end{proof}

\subsection{Archimedean estimates for \texorpdfstring{$\boldsymbol{d>2}$}{$d>2$}}\label{sec4.3}

In this section, we assume that $K$ is a number field, $a,b\in K$,
$v\in M_{K}^{\infty}$, and $d>2$. Because $v$ is fixed, we write
$|\cdot|$, $M_{a}$, $\mu_{a}$, and $g_{a}$ for $|\cdot|_{v}$,
$M_{a,v}$, $\mu_{a,v}$, and $g_{a,v}$.

We will prove the counterparts of Propositions~\ref{prop4.4},~\ref{prop4.5},
and~\ref{prop4.6} for $d>2$, but this time we only need to consider
the roots of $f_{T}^{2}(a)=f_{T}(a)$. The reason is as follows: In
Equations (\ref{eq4.1}) and (\ref{eq4.2}), what we really need
is
\[
|a^{d}-b^{d}|\geq\frac{c_{1}(d)}{\max\{|a|,|b|\}^{c_{2}(d)}}
\]
for some $c_{1}(d),c_{2}(d)>0$. When we consider the roots of $f_{T}^{2}(a)=f_{T}(a)$,
we will get $c_{2}(d)=d-2$. This is already good enough for $d>2$,
so we no longer need to consider the roots of $f_{T}^{3}(a)=f_{T}(a)$.
For $0\leq i\leq d-1$, let
\[
\alpha_{i}(a)=-a^{d}+\zeta_{d}^{i}a
\]
be the roots of $f_{T}^{2}(a)=f_{T}(a)$.

\begin{prop}\label{prop4.7}
  Assume that $|a|\geq6$. Then we have
\begin{enumerate}
\item $M_{a}\subseteq\cup_{i=0}^{d-1}D(\alpha_{i}(a),12/|a|^{d-2})$; 
\item if $t\notin\cup_{i=0}^{d-1}D(\alpha_{i}(a),12/|a|^{d-2})$, then
\[
g_{a}(t)\geq\frac{1}{d}\log(13|a|); 
\]
\item $\mu_{a}(D(\alpha_{i}(a),12/|a|^{d-2}))=1/d$ for $0\leq i\leq d-1$.
\end{enumerate}
\end{prop}

\begin{proof}
For simplicity, let $\alpha_{i}=\alpha_{i}(a)$ and $D_{i}=D(\alpha_{i},12/|a|^{d-2})$.
If $t\in\partial D_{0}$, then
\begin{align*}
|t-\alpha_{0}| & =12/|a|^{d-2},\\
|t-\alpha_{i}| & \geq|\alpha_{0}-\alpha_{i}|-|t-\alpha_{0}|=\left|1-\zeta_{d}^{i}\right||a|-12/|a|^{d-2}\geq(1-1/d)\left|1-\zeta_{d}^{i}\right||a|,
\end{align*}
where the last inequality follows from
\[
\left(\frac{12d}{\left|1-\zeta_{d}^{i}\right|}\right)^{\frac{1}{d-1}}\leq\left(\frac{12d^{2}}{d\left|1-\zeta_{d}\right|}\right)^{\frac{1}{d-1}}\leq\left(\frac{4d^{2}}{\sqrt{3}}\right)^{\frac{1}{d-1}}\leq2\cdot3^{3/4}\leq|a|.
\]
Then we have
\[
|f_{t}^{2}(a)-f_{t}(a)|\geq\frac{12}{|a|^{d-2}}\prod_{i=1}^{d-1}\left(\left(1-\frac{1}{d}\right)|1-\zeta_{d}^{i}|\cdot|a|\right)=12d\left(1-\frac{1}{d}\right)^{d-1}|a|\geq16|a|.
\]
The rest of the proof is similar to that of Proposition~\ref{prop4.4}.
\end{proof}

\begin{prop}\label{prop4.8}
For any $a,b\in K$, we have
\[
\int g_{a}d\mu_{b}\geq\frac{1}{d^{3}}\log^{+}\left|a^{d}-b^{d}\right|-\frac{1}{d^{3}}\log\left(2\cdot9^{d}\right).
\]
Moreover, if\, $\max\{|a|,|b|\}\geq9$ and
\begin{equation}
\left|a^{d}-b^{d}\right|\geq\frac{25}{\max\{|a|,|b|\}^{d-2}},\label{eq4.2}
\end{equation}
then
\[
\int g_{a}d\mu_{b}\geq\frac{1}{d^{2}}\log\max\{|a|,|b|\}.
\]
\end{prop}

\begin{proof}
Without loss of generality, we assume that $|a|\geq|b|$. We prove
the second assertion in parts~\eqref{4.8-part1},~\eqref{4.8-part2}, and we prove the first assertion
in part~\eqref{4.8-part3}.

\begin{enumerate}[wide]
\item\label{4.8-part1}
  \fbox{$|a|\geq9$ and $|b|\leq6$} If $t\in M_{b}$, then by
Proposition~\ref{prop3.4},
\[
|t|\leq3\max\{|b|,4\}^{d}\leq3\cdot6^{d}\leq9\left(9^{d-1}-2\right)\leq|a|\left(|a|^{d-1}-2\right)=|a|^{d}-2|a|
\]
and
\[
|f_{t}(a)|\geq|a|^{d}-|t|\geq2|a|.
\]
By an argument similar to that of the proof of Proposition~\ref{prop4.6},
we have
\[
\int g_{a}d\mu_{b}\geq\log|a|.
\]

\item\label{4.8-part2}
\fbox{$|a|\geq9$, $|b|\geq6$, and $|a^{d}-b^{d}|\geq25/|a|^{d-2}$}

\begin{claim*} Some $D(\alpha_{i}(b),12/|b|^{d-2})$ is disjoint
from $\cup_{i=0}^{d-1}D(\alpha_{i}(a),12/|a|^{d-2})$.
\end{claim*}

Suppose the claim is false; then for each $i$, there exists a $k_{i}$ such that
\[
|\alpha_{k_{i}}(a)-\alpha_{i}(b)|<12/|a|^{d-2}+12/|b|^{d-2}.
\]
If $k_{i}=k_{j}$ for some $i\neq j$, then
\begin{align*}
\left|\zeta_{d}^{i}-\zeta_{d}^{j}\right||b| & =|\alpha_{i}(b)-\alpha_{j}(b)|\leq|\alpha_{k_{i}}(a)-\alpha_{i}(b)|+|\alpha_{k_{j}}(a)-\alpha_{j}(b)|\\
 & <24/|a|^{d-2}+24/|b|^{d-2}\leq48/|b|^{d-2}.
\end{align*}
This gives a contradiction because
\[
\left(\frac{48}{\left|\zeta_{d}^{i}-\zeta_{d}^{j}\right|}\right)^{\frac{1}{d-1}}\leq\left(\frac{48d}{d\left|1-\zeta_{d}\right|}\right)^{\frac{1}{d-1}}\leq\left(\frac{16d}{\sqrt{3}}\right)^{\frac{1}{d-1}}\leq4\cdot3^{1/4}\leq|b|.
\]
If $k_{i}\neq k_{j}$ for any $i\neq j$, then
\[
d|a^{d}-b^{d}|=\left|\sum_{i=0}^{d-1}\alpha_{i}(a)-\sum_{i=0}^{d-1}\alpha_{i}(b)\right|\leq\sum_{i=0}^{d-1}\left|\alpha_{k_{i}}(a)-\alpha_{i}(b)\right|<d\left(12/|a|^{d-2}+12/|b|^{d-2}\right).
\]
If $|a^{d}-b^{d}|\geq4$, then
\[
4\leq\left|a^{d}-b^{d}\right|<12/|a|^{d-2}+12/|b|^{d-2}\leq24/|b|^{d-2}\leq4,
\]
which is false. If $25/|a|^{d-2}\leq|a^{d}-b^{d}|\leq4$, then
\[
\frac{|a|^{d-2}}{|b|^{d-2}}\leq\frac{|a|^{d}}{|b|^{d}}\leq\frac{|a^{d}-b^{d}|}{|b|^{d}}+1\leq\frac{4}{6^{d}}+1\leq\frac{13}{12}
\]
and
\[
25/|a|^{d-2}\leq|a^{d}-b^{d}|<12/|a|^{d-2}+12/|b|^{d-2}\leq25/|a|^{d-2},
\]
which is also false. Therefore, the claim is proved and, by Proposition
\ref{prop4.7},
\[
\int g_{a}d\mu_{b}\geq\frac{1}{d}\cdot\frac{1}{d}\log(13|a|)=\frac{1}{d^{2}}\log(13|a|).
\]

\item\label{4.8-part3}
Now we prove the first assertion. If $|a|\leq9$ or $|a^{d}-b^{d}|\leq3$,
then
\[
\log^{+}\left|a^{d}-b^{d}\right|\leq\log\max\left\{2|a|^{d},3\right\}\leq\log\left(2\cdot9^{d}\right)
\]
and
\[
\int g_{a}d\mu_{b}\geq0\geq\frac{1}{d^{3}}\log^{+}\left|a^{d}-b^{d}\right|-\frac{1}{d^{3}}\log\left(2\cdot9^{d}\right).
\]
If $|a|\geq9$ and $|a^{d}-b^{d}|\geq3\geq25/|a|^{d-2}$, then by
parts~\eqref{4.8-part1} and~\eqref{4.8-part2},
\[
\int g_{a}d\mu_{b}\geq\frac{1}{d^{2}}\log(13|a|)\geq\frac{1}{d^{3}}\log\left(2|a|^{d}\right)\geq\frac{1}{d^{3}}\log^{+}\left|a^{d}-b^{d}\right|.\qedhere
\]
\end{enumerate}
\end{proof}

\subsection{Non-Archimedean estimates}\label{sec4.4}

In this section, we assume that $K$ is a number field, $a,b\in K$,
and $v\in M_{K}^{0}$. Because $v$ is fixed, we write $|\cdot|$,
$\delta$, $M_{a}$, $\mu_{a}$, and $g_{a}$ for $|\cdot|_{v}$,
$\delta_{v}$, $M_{a,v}$, $\mu_{a,v}$, and $g_{a,v}$.

As in \cite[Sections 5.1 and 6.1]{DKY2}, we first study the structure
of $M_{a}$. More precisely, we show that when~$a$ is large enough,
$M_{a}$ can be described with respect to the roots of $f_{T}^{n}(a)$.

\begin{prop}\label{prop4.9}
Assume that $|a|>|d|^{-2/(d-1)}$. Fix a $t\in\C_{v}$
such that $|f_{t}^{n}(a)|\leq|a|$ for some $n\geq1$. Let $s_{1},\dots,s_{d^{n-1}}$
be the roots of $f_{T}^{n}(a)$ such that
\[
|t-s_{1}|\leq\dots\leq|t-s_{d^{n-1}}|.
\]
Then we have
\[
|t-s_{1}|=\frac{\left|f_{t}^{n}(a)\right|}{\left(|d|\cdot|a|^{d-1}\right)^{n-1}}\leq\frac{|a|}{\left(|d|\cdot|a|^{d-1}\right)^{n-1}}<|t-s_{2}|\leq\frac{|a|}{\left(|d|\cdot|a|^{d-1}\right)^{n-2}}<|t-s_{d+1}|.
\]
\end{prop}

\begin{proof}
Let $A_{n}=|a|/(|d|\cdot|a|^{d-1})^{n-1}$ for any $n\geq1$. By Proposition
\ref{prop3.10}, $|f_{t}^{n}(a)|\leq|a|$ implies $|f_{t}^{k}(a)|=|a|$
for any $1\leq k\leq n-1$. Thus we can prove the assertion by induction.
Because $|t+a^{d}|=|f_{t}(a)|$, the statement is true for $n=1$.
Assume the statement is true for some $n\geq1$, and fix a $t\in\C_{v}$
such that $|f_{t}^{n+1}(a)|\leq|a|$. Let $\alpha_{1},\dots,\alpha_{d^{n-1}}$
be the roots of $f_{T}^{n}(a)$ with $|t-\alpha_{i}|$ increasing,
and let $\beta_{1},\dots,\beta_{d^{n}}$ be the roots of $f_{T}^{n+1}(a)$
with $|t-\beta_{i}|$ increasing. Then
\[
\prod_{i=1}^{d^{n-1}}(t-T-\alpha_{i})^{d}+t-T=f_{t-T}^{n}(a)^{d}+t-T=f_{t-T}^{n+1}(a)=\prod_{i=1}^{d^{n}}(t-T-\beta_{i}).
\]
Let
\[
\prod_{i=1}^{d^{n-1}}(t-T-\alpha_{i})^{d}=\sum_{i=0}^{d^{n}}a_{i}T^{i}\text{ and }\prod_{i=1}^{d^{n}}(t-T-\beta_{i})=\sum_{i=0}^{d^{n}}b_{i}T^{i};
\]
then by the assumption $|a|>|d|^{-2/(d-1)}$ and the induction hypothesis
$|t-\alpha_{1}|=A_{n}<|t-\alpha_{2}|$,
\begin{align*}
|b_{0}| & =\left|\prod_{i=1}^{d^{n}}(t-\beta_{i})\right|=\left|f_{t}^{n+1}(a)\right|,\\
|b_{1}| & =|a_{1}-1|=\left|\sum_{j=1}^{d^{n-1}}\frac{d\prod_{i=1}^{d^{n-1}}(t-\alpha_{i})^{d}}{t-\alpha_{j}}-1\right|=\frac{|d|\left|f_{t}^{n}(a)\right|^{d}}{A_{n}}=\frac{|d|\cdot|a|^{d}}{A_{n}},\\
|b_{i}| & =|a_{i}|\leq\frac{\left|f_{t}^{n}(a)\right|^{d}}{A_{n}^{i}}=\frac{|a|^{d}}{A_{n}^{i}}\quad\text{for any }2\leq i\leq d-1,\\
|b_{d}| & =|a_{d}|=\frac{\left|f_{t}^{n}(a)\right|^{d}}{A_{n}^{d}}=\frac{|a|^{d}}{A_{n}^{d}}.
\end{align*}
Now we consider the Newton polygon of $f_{t-T}^{n+1}(a)$. By the
induction hypothesis,
\[
\left(\frac{|b_{i}|}{|b_{d}|}\right)^{\frac{1}{d-i}}\leq A_{n}<\left(\frac{|a_{d}|}{|a_{j}|}\right)^{\frac{1}{j-d}}=\left(\frac{|b_{d}|}{|b_{j}|}\right)^{\frac{1}{j-d}}
\]
for any $i<d<j$, so $(d,-\log|b_{d}|)$ is a vertex of the Newton
polygon. By the assumptions,
\[
\frac{|b_{0}|}{|b_{1}|}=\frac{\left|f_{t}^{n+1}(a)\right|A_{n}}{|d|\cdot|a|^{d}}\leq\frac{A_{n}}{|d|\cdot|a|^{d-1}}<|d|A_{n}\leq\left(\frac{|b_{1}|}{|b_{i}|}\right)^{\frac{1}{i-1}}
\]
for any $2\leq i\leq d$, so $(1,-\log|b_{1}|)$ is also a vertex
of the Newton polygon. Therefore,
\[
|t-\beta_{1}|=\frac{\left|f_{t}^{n+1}(a)\right|A_{n}}{|d|\cdot|a|^{d}}\leq\frac{A_{n}}{|d|\cdot|a|^{d-1}}<|t-\beta_{2}|\leq A_{n}<|t-\beta_{d+1}|.
\]
This completes the inductive step and hence the proof.
\end{proof}

In particular, Proposition~\ref{prop4.9} implies the non-Archimedean
versions of Propositions~\ref{prop4.4},~\ref{prop4.5}, and~\ref{prop4.7}
as follows.

\begin{prop}\label{prop4.10}
  Assume that $|a|>|d|^{-2/(d-1)}$. For any $n\geq1$,
let $S_{n}$ be the set of all roots of $f_{T}^{n}(a)$. Then  
\begin{enumerate}
\item\label{p4.10-1} all roots of $f_{T}^{n}(a)$ are simple; moreover, for any $s_{1},s_{2}\in S_{n}$,
we have
\[
|s_{1}-s_{2}|>\frac{|a|}{\left(|d|\cdot|a|^{d-1}\right)^{n-1}};
\]
\item\label{p4.10-2}  $M_{a}\subseteq\C_{v}$;
\item\label{p4.10-3} for any $s\in S_{n}$, we have
\[
\mu_{a}\left(\mathcal{D}\left(s,\frac{|a|}{\left(|d|\cdot|a|^{d-1}\right)^{n-1}}\right)\right)=\frac{1}{d^{n-1}}.
\]
\end{enumerate}
\end{prop}

\begin{proof}
For any $n\geq1$, let $A_{n}=|a|/(|d|\cdot|a|^{d-1})^{n-1}$. For any
$s\in S_{n}$, let
\[
T_{n,s}=\{t\in\C_{v}:|s-t|=A_{n}\}.
\]

\eqref{p4.10-1}~
 Fix any $s\in S_{n}$, and let $s_{1},\dots,s_{d^{n-1}}\in S_{n}$
be such that $|s-s_{i}|$ is increasing. By Proposition~\ref{prop4.9},
we have $|s-s_{1}|=0<A_{n}<|s-s_{2}|$, so $s$ is simple.

\begin{claim}\label{claim1}
For any $s\in S_{n}$, we have $|S_{n+1}\cap T_{n,s}|=d$.
\end{claim}

By Propositions~\ref{prop3.10} and~\ref{prop4.9}, we have
$S_{n+1}\subseteq\cup_{s\in S_{n}}T_{n,s}$.  Since
$|S_{n+1}|=d|S_{n}|$, it suffices to show that $|S_{n+1}\cap
T_{n,s}|\leq d$ for any $s\in S_{n}$. Suppose
$\alpha_{1},\dots,\alpha_{d+1}\in S_{n+1}\cap T_{n,s}$; then for any
$2\leq i\leq d+1$, we have
\[
|\alpha_{1}-\alpha_{i}|\leq\max\{|s-\alpha_{1}|,|s-\alpha_{i}|\}=A_{n},
\]
which contradicts Proposition~\ref{prop4.9}.

\begin{claim}\label{claim2}
If $t\in T_{n,s}$ for some $s\in S_{n}$, then
$|f_{t}^{n}(a)|=|a|$.
\end{claim}

By Claim~\ref{claim1}, there exists an $\alpha\in S_{n+1}\cap T_{n,s}$. By Proposition~\ref{prop4.9}, for any $s'\in S_{n}\backslash\{s\}$, we have
\begin{align*}
|s-t| & =A_{n}=|s-\alpha|,\\
|s'-t| & =\max\{|s'-\alpha|,|s-\alpha|,|s-t|\}=|s'-\alpha|.
\end{align*}
By Proposition~\ref{prop3.10}, we have
\[
|f_{t}^{n}(a)|=|s-t|\prod_{s'\in S_{n}\backslash\{s\}}|s'-t|=|s-\alpha|\prod_{s'\in S_{n}\backslash\{s\}}|s'-\alpha|=|f_{\alpha}^{n}(a)|=|a|.
\]

\eqref{p4.10-2}~ Let $T_{n}=\{t\in\C_{v}:|f_{t}^{n}(a)|=|a|\}$; then by
Proposition~\ref{prop3.9},
\[
M_{a}\cap\C_{v}=\cap_{n=1}^{\infty}T_{n}\subseteq\cap_{n=1}^{\infty}\overline{T_{n}}.
\]
By Proposition~\ref{prop4.9} and Claim~\ref{claim2},
\begin{align*}
T_{n} & =\{t\in\C_{v}:|s-t|=A_{n}\text{ for some }s\in S_{n}\},\\
\overline{T_{n}}& =\left\{t\in\A_{v}^{1,\an}:\delta(s,t)=A_{n}\text{ for some }s\in S_{n}\right\}.
\end{align*}
If $t\in\cap_{n=1}^{\infty}\overline{T_{n}}$, then $\diam(t)\leq A_{n}$
for any $n\geq1$. Since $A_{n}\to0$ as $n\to\infty$, we have $t\in\C_{v}$.
Therefore, $\cap_{n=1}^{\infty}T_{n}=\cap_{n=1}^{\infty}\overline{T_{n}}$
and $M_{a}\cap\C_{v}=\overline{M_{a}\cap\C_{v}}=M_{a}$.

\eqref{p4.10-3}~
 By Claim~\ref{claim1} and induction, for any $s\in S_{n}$ and any $k\geq n+1$,
we have $|S_{k}\cap T_{n,s}|=d^{k-n}$. Then the conclusion follows
from Proposition~\ref{prop4.2}.
\end{proof}

Propositions~\ref{prop4.11} and~\ref{prop4.12} can be seen as the
non-Archimedean versions of Propositions~\ref{prop4.6} and~\ref{prop4.8}.

\begin{prop}\label{prop4.11}
Assume that $\max\{|a|,|b|\}>|d|^{-2/(d-1)}$. If
\[
|a^{d}-b^{d}|>\frac{|d|^{-2(n-1)}}{\max\{|a|,|b|\}^{(d-1)(n-1)-1}} 
\]
for some $n\geq1$, then
\[
\int g_{a}d\mu_{b}\geq\frac{1}{d^{2n-2}}\log\max\{|a|,|b|\}.
\]
\end{prop}

\begin{proof}
If $|a|\neq|b|$, then by Propositions~\ref{prop3.8} and~\ref{prop3.9},
\[
\int g_{a}d\mu_{b}=d\log\max\{|a|,|b|,1\}\geq\frac{1}{d^{2n-2}}\log\max\{|a|,|b|\}.
\]
From now on we assume that $|a|=|b|$. For any $n\geq1$, let $A_{n}=|a|/(|d|\cdot|a|^{d-1})^{n-1}$,
and let $\alpha_{1}(a),\dots,\alpha_{d^{n-1}}(a)$ be the roots of
$f_{T}^{n}(a)$. For any $1\leq i\leq d^{n-1}$, let $D_{i}(a)=\overline{D}(\alpha_{i}(a),A_{n})$.
Define $\alpha_{i}(b)$ and $D_{i}(b)$ similarly.

\begin{claim*}
  Some $D_{i}(b)$ is disjoint from $\cup_{i=1}^{d^{n-1}}D_{i}(a)$.
\end{claim*}

Suppose the claim is false; then for each $i$, there exists a $k_{i}$ such that
\[
|\alpha_{k_{i}}(a)-\alpha_{i}(b)|\leq A_{n}.
\]
If $k_{i}=k_{j}$ for some $i\neq j$, then
\[
|\alpha_{i}(b)-\alpha_{j}(b)|\leq\max\left\{\left|\alpha_{k_{i}}(a)-\alpha_{i}(b)\right|,\left|\alpha_{k_{j}}(a)-\alpha_{j}(b)\right|\right\}\leq A_{n},
\]
which contradicts Proposition~\ref{prop4.10}. If $k_{i}\neq k_{j}$
for any $i\neq j$, then
\[
|d|^{n-1}\left|a^{d}-b^{d}\right|=\left|\sum_{i=1}^{d^{n-1}}\alpha_{i}(a)-\sum_{i=1}^{d^{n-1}}\alpha_{i}(b)\right|\leq\max\{|\alpha_{k_{i}}(a)-\alpha_{i}(b)|\}_{i=1}^{d^{n-1}}\leq A_{n},
\]
which contradicts the assumption
\[
\left|a^{d}-b^{d}\right|>\frac{|d|^{-2(n-1)}}{|a|^{(d-1)(n-1)-1}}=\frac{|a|}{\left(|d|^{2}|a|^{d-1}\right)^{n-1}}=\frac{A_{n}}{|d|^{n-1}}.
\]
Let $D_{i}(b)$ be disjoint from $\cup_{i=1}^{d^{n-1}}D_{i}(a)$ and
$t\in D_{i}(b)$. By Proposition~\ref{prop3.10},
\[
|t|=\max\{|t-\alpha_{i}(b)|,|\alpha_{i}(b)|\}=|a|^{d}.
\]
Since $t\notin\cup_{i=1}^{d^{n-1}}D_{i}(a)$, by Proposition~\ref{prop4.9},
we have $|f_{t}^{n}(a)|>|a|$. By induction, $|f_{t}^{n+k}(a)|=|f_{t}^{n}(a)|^{d^{k}}$
for any $k\geq1$. Therefore, by Theorem~\ref{thm2.1},
\[
g_{a}(t)=\lim_{k\to\infty}\frac{1}{d^{n+k-1}}\log^{+}\left|f_{t}^{n+k}(a)\right|=\frac{1}{d^{n-1}}\log\left|f_{t}^{n}(a)\right|\geq\frac{1}{d^{n-1}}\log|a|.
\]
By Proposition~\ref{prop4.10},
\[
\int g_{a}d\mu_{b}\geq\frac{1}{d^{n-1}}\cdot\frac{1}{d^{n-1}}\log|a|=\frac{1}{d^{2n-2}}\log|a|.\qedhere
\]
\end{proof}

\begin{prop}\label{prop4.12}
For any $a,b\in K$, we have
\[
\int g_{a}d\mu_{b}\geq\frac{1}{d^{3}}\log^{+}\left|a^{d}-b^{d}\right|-\frac{1}{d^{3}}\log|d|^{-2d/(d-1)}.
\]
\end{prop}

\begin{proof}
If $\max\{|a|,|b|\}\leq|d|^{-2/(d-1)}$ or $|a^{d}-b^{d}|\leq|d|^{-2/(d-1)}$,
then
\[
\log^{+}\left|a^{d}-b^{d}\right|\leq\log\max\left\{|a|^{d},|b|^{d},|d|^{-2/(d-1)}\right\}\leq\log|d|^{-2d/(d-1)}
\]
and
\[
\int g_{a}d\mu_{b}\geq0\geq\frac{1}{d^{3}}\log^{+}\left|a^{d}-b^{d}\right|-\frac{1}{d^{3}}\log|d|^{-2d/(d-1)}.
\]
If $\max\{|a|,|b|\}>|d|^{-2/(d-1)}$ and
\[
\left|a^{d}-b^{d}\right|>|d|^{-2/(d-1)}=\frac{|d|^{-2}}{|d|^{-2(d-2)/(d-1)}}\geq\frac{|d|^{-2}}{\max\{|a|,|b|\}^{d-2}},
\]
then by Proposition~\ref{prop4.11},
\[
\int g_{a}d\mu_{b}\geq\frac{1}{d^{2}}\log\max\{|a|,|b|\}=\frac{1}{d^{3}}\log\max\left\{|a|^{d},|b|^{d}\right\}\geq\frac{1}{d^{3}}\log^{+}\left|a^{d}-b^{d}\right|.\qedhere
\]
\end{proof}

\subsection{Proofs of Theorems \texorpdfstring{\ref{thm1.3}}{1.3}
and \texorpdfstring{\ref{thm1.4}}{1.4}}\label{sec4.5}

Now we are ready to give the proofs of Theorems~\ref{thm1.3} and
\ref{thm1.4}. Given Propositions~\ref{prop4.6},~\ref{prop4.8},
\ref{prop4.11}, and~\ref{prop4.12}, our proofs are almost identical
to those given in \cite[Sections 7 and 8]{DKY2}.

\begin{proof}[Proof of Theorem~\ref{thm1.3}]
Let $K$ be a number field such
that $a,b\in K$. Define
\[
(p_{v},q_{v},m,n)=\begin{cases}
(50,11,1,4) & \text{if }d=2\text{ and }v\in M_{K}^{\infty},\\
\left(|2|_{v}^{-2},|2|_{v}^{-4},1,4\right) & \text{if }d=2\text{ and }v\in M_{K}^{0},\\
(9,25,d-2,2) & \text{if }d>2\text{ and }v\in M_{K}^{\infty},\\
\left(|d|_{v}^{-2/(d-1)},|d|_{v}^{-2},d-2,2\right) & \text{if }d>2\text{ and }v\in M_{K}^{0},
\end{cases}
\]
and $r_{v}=\max\{|a|_{v},|b|_{v},p_{v}\}$ for any $v\in M_{K}$.
Following \cite[Section 8.2]{DKY2}, we define
\begin{align*}
M_{\help} & =\left\{ v\in M_{K}:\max\{|a|_{v},|b|_{v}\}>p_{v},\ \left|a^{d}-b^{d}\right|_{v}>\frac{q_{v}}{\max\{|a|_{v},|b|_{v}\}^{m}}\right\} ,\\
M_{\close} & =\left\{ v\in M_{K}:\max\{|a|_{v},|b|_{v}\}>p_{v},\ \left|a^{d}-b^{d}\right|_{v}\leq\frac{q_{v}}{\max\{|a|_{v},|b|_{v}\}^{m}}\right\} ,\\
M_{\bounded} & =\{v\in M_{K}:\max\{|a|_{v},|b|_{v}\}\leq p_{v}\}.
\end{align*}
Then
\allowdisplaybreaks
\begin{align*}
0 & =\sum_{v\in M_{K}}n_{v}\log\left|a^{d}-b^{d}\right|_{v}\\
 & \leq\sum_{v\in M_{\close}}n_{v}\log\frac{q_{v}}{r_{v}^{m}}+\sum_{v\in M_{K}^{\infty}\backslash M_{\close}}n_{v}\log\left(2r_{v}^{d}\right)+\sum_{v\in M_{K}^{0}\backslash M_{\close}}n_{v}\log r_{v}^{d}\\
 & \leq\sum_{v\in M_{\close}}n_{v}\log\frac{q_{v}}{r_{v}^{m}}+\sum_{v\in M_{K}\backslash M_{\close}}n_{v}\log\left(q_{v}r_{v}^{d}\right)\\
 & =\sum_{v\in M_{K}}n_{v}\log q_{v}-\sum_{v\in M_{\close}}n_{v}\log r_{v}^{m}+\sum_{v\in M_{K}\backslash M_{\close}}n_{v}\log r_{v}^{d}\\
 & =\sum_{v\in M_{K}}n_{v}\log q_{v}-\sum_{v\in M_{K}}n_{v}\log r_{v}^{m}+\sum_{v\in M_{\help}}n_{v}\log r_{v}^{d+m}+\sum_{v\in M_{\bounded}}n_{v}\log r_{v}^{d+m}\\
 & \leq\sum_{v\in M_{K}}n_{v}\log q_{v}^{d+m}-\sum_{v\in M_{K}}n_{v}\log r_{v}^{m}+\sum_{v\in M_{\help}}n_{v}\log r_{v}^{d+m}+\sum_{v\in M_{K}}n_{v}\log p_{v}^{d+m}\\
 & \leq(d+m)\sum_{v\in M_{\help}}n_{v}\log r_{v}-mh(a,b)+(d+m)\sum_{v\in M_{K}}n_{v}\log(p_{v}q_{v}).
\end{align*}
By Propositions~\ref{prop4.6},~\ref{prop4.8}, and~\ref{prop4.11},
\begin{align*}
\left\langle \boldsymbol{\mu}_{a},\boldsymbol{\mu}_{b}\right\rangle  & \geq\sum_{v\in M_{\help}}n_{v}\int g_{a,v}d\mu_{b,v}\geq\frac{1}{d^{n}}\sum_{v\in M_{\help}}n_{v}\log r_{v}\\
 & \geq\frac{1}{d^{n}}\left(\frac{m}{d+m}h(a,b)-\sum_{v\in M_{K}}n_{v}\log(p_{v}q_{v})\right)\\
 & =\begin{cases}
\frac{1}{48}h(a,b)-\frac{1}{16}\log35200 & \text{if }d=2,\\
\frac{d-2}{2d^{2}(d-1)}h(a,b)-\frac{1}{d^{2}}\log(225d^{2d/(d-1)}) & \text{if }d>2,
\end{cases}\\
 & \geq\frac{1}{12d^{2}}h(a,b)-1.\qedhere
\end{align*}
\end{proof}

The proof of Theorem~\ref{thm1.4} relies on the following two results.
Theorem~\ref{thm4.13} is adapted from \cite[Theorem 7.1]{DKY2},
and Proposition~\ref{prop4.14} is adapted from a continuity argument
used in the proof of \cite[Theorem 1.6]{DKY2}.

\begin{thm}
\label{thm4.13}For any $a,b\in\bar{\Q}$, we have
\[
\left\langle \boldsymbol{\mu}_{a},\boldsymbol{\mu}_{b}\right\rangle \geq\frac{1}{d^{3}}h\left(a^{d}-b^{d}\right)-2,
\]
where $h$ is the logarithmic Weil height on $\bar{\Q}$.
\end{thm}

\begin{proof}
Let $K$ be a number field such that $a,b\in K$. Define
\[
r_{v}=\begin{cases}
5000 & \text{if }d=2\text{ and }v\in M_{K}^{\infty},\\
|2|_{v}^{-4} & \text{if }d=2\text{ and }v\in M_{K}^{0},\\
2\cdot9^{d} & \text{if }d>2\text{ and }v\in M_{K}^{\infty},\\
|d|_{v}^{-2d/(d-1)} & \text{if }d>2\text{ and }v\in M_{K}^{0}.
\end{cases}
\]
By Propositions~\ref{prop4.6},~\ref{prop4.8}, and~\ref{prop4.12},
\allowdisplaybreaks
\begin{align*}
\left\langle \boldsymbol{\mu}_{a},\boldsymbol{\mu}_{b}\right\rangle  & =\sum_{v\in M_{K}}n_{v}\int g_{a,v}d\mu_{b,v}\geq\sum_{v\in M_{K}}n_{v}\left(\frac{1}{d^{3}}\log^{+}\left|a^{d}-b^{d}\right|_{v}-\frac{1}{d^{3}}\log r_{v}\right)\\
 & =\frac{1}{d^{3}}h\left(a^{d}-b^{d}\right)-\frac{1}{d^{3}}\sum_{v\in M_{K}}n_{v}\log r_{v}\\
& =\frac{1}{d^{3}}h\left(a^{d}-b^{d}\right)-
\begin{cases}
\frac{1}{8}\log80000 & \text{if }d=2,\\
\frac{1}{d^{3}}\log\left(2\cdot9^{d}d^{2d/(d-1)}\right) & \text{if }d>2,
\end{cases}\\
 & \geq\frac{1}{d^{3}}h\left(a^{d}-b^{d}\right)-2.\qedhere
\end{align*}
\end{proof}

\begin{prop}\label{prop4.14}
The complex Arakelov--Zhang pairing $\int g_{a}d\mu_{b}$
is a continuous function of $(a,b)\in\C^{2}$.
\end{prop}

\begin{proof}
Let $a_{n}\to a$ and $b_{n}\to b$ as $n\to\infty$; then by Proposition
\ref{prop3.4}, there exists an $r>0$ such that for any $n\geq1$,
\[
\{a,b,a_{n},b_{n}\},M_{a},M_{b},M_{a_{n}},M_{b_{n}}\subseteq\overline{D}(0,r).
\]
By \cite[Proposition 1.2]{BH}, $g_{a}(t)$ is a continuous function
of $(a,t)\in\C^{2}$, so $g_{a}(t)$ is uniformly continuous
on the compact set $\overline{D}(0,r)\times\overline{D}(0,r)$. Therefore,
as $n\to\infty$,
\[
\int g_{a}d\mu_{b}-\int g_{a_{n}}d\mu_{b_{n}}=\int(g_{a}-g_{a_{n}})d\mu_{b}+\int(g_{b}-g_{b_{n}})d\mu_{a_{n}}\longrightarrow 0.\qedhere
\]
\end{proof}

\begin{proof}[Proof of Theorem~\ref{thm1.4}]
Let $K$ be a number field such that $a,b\in K$. Define
\[
(p,q,m,n)=\begin{cases}
(50,11,1,4) & \text{if }d=2,\\
(9,25,d-2,2) & \text{if }d>2.
\end{cases}
\]
Let $N$ be a large number to be determined later, and let
\begin{align*}
M_{1} & =\left\{v\in M_{K}^{\infty}:\max\{|a|_{v},|b|_{v}\}\leq N,\ \left|a^{d}-b^{d}\right|_{v}\geq q/N^{m}\right\},\\
M_{2} & =\left\{v\in M_{K}^{\infty}:\max\{|a|_{v},|b|_{v}\}>N,\ \left|a^{d}-b^{d}\right|_{v}\geq q/N^{m}\right\},\\
M_{3} & =\left\{v\in M_{K}^{\infty}:\left|a^{d}-b^{d}\right|_{v}<q/N^{m}\right\}.
\end{align*}
Since $\sum_{v\in M_{K}^{\infty}}n_{v}=1$, we have $\sum_{v\in M_{i}}n_{v}\geq1/3$
for some $1\leq i\leq3$.

\begin{enumerate}[wide]
\item Assume that $\sum_{v\in M_{1}}n_{v}\geq1/3$. Let
\[
R_{N}=\left\{(c_{1},c_{2})\in\C^{2}:\max\{|c_{1}|,|c_{2}|\}\leq N,\ \left|c_{1}^{d}-c_{2}^{d}\right|\geq q/N^{m}\right\}.
\]
Since $R_{N}$ is compact, by \cite[Lemma 3.4]{BD1} and Proposition
\ref{prop4.14}, we have 
\[
r_{N}=\min_{(c_{1},c_{2})\in R_{N}}\int g_{c_{1}}d\mu_{c_{2}}>0.
\]
Therefore,
\[
\left\langle \boldsymbol{\mu}_{a},\boldsymbol{\mu}_{b}\right\rangle \geq\sum_{v\in M_{1}}n_{v}\int g_{a,v}d\mu_{b,v}\geq\frac{1}{3}r_{N}.
\]

\item Assume that $\sum_{v\in M_{2}}n_{v}\geq1/3$. If $v\in M_{2}$,
then
\[
|a^{d}-b^{d}|_{v}\geq\frac{q}{N^{m}}\geq\frac{q}{\max\{|a|_{v},|b|_{v}\}^{m}}.
\]
If $N\geq p$, then by Propositions~\ref{prop4.6} and~\ref{prop4.8},
\[
\int g_{a,v}d\mu_{b,v}\geq\frac{1}{d^{n}}\log\max\{|a|_{v},|b|_{v}\}\geq\frac{1}{d^{n}}\log N.
\]
Therefore,
\[
\left\langle \boldsymbol{\mu}_{a},\boldsymbol{\mu}_{b}\right\rangle \geq\sum_{v\in M_{2}}n_{v}\int g_{a,v}d\mu_{b,v}\geq\frac{1}{3d^{n}}\log N.
\]

\item Assume that $\sum_{v\in M_{3}}n_{v}\geq1/3$. By Theorem~\ref{thm4.13},
\allowdisplaybreaks
  \begin{align*}
\left\langle \boldsymbol{\mu}_{a},\boldsymbol{\mu}_{b}\right\rangle  & \geq\frac{1}{d^{3}}h\left(a^{d}-b^{d}\right)-2=\frac{1}{d^{3}}\sum_{v\in M_{K}}n_{v}\log^{+}\left|a^{d}-b^{d}\right|_{v}-2\\
 & \geq\frac{1}{d^{3}}\sum_{v\in M_{K}\backslash M_{3}}n_{v}\log^{+}\left|a^{d}-b^{d}\right|_{v}-2\geq\frac{1}{d^{3}}\sum_{v\in M_{K}\backslash M_{3}}n_{v}\log\left|a^{d}-b^{d}\right|_{v}-2\\
 & =\frac{1}{d^{3}}\sum_{v\in M_{3}}n_{v}\log\frac{1}{\left|a^{d}-b^{d}\right|_{v}}-2\geq\frac{1}{3d^{3}}\log\frac{N^{m}}{q}-2,
\end{align*}
  which is positive when $N>q^{1/m}e^{6d^{3}/m}$.\hfill\qedhere
  \end{enumerate}
\end{proof}

\section{Proof of Theorem \texorpdfstring{\ref{thm1.1}}{1.1}}\label{sec5}

Finally, we deduce Theorem~\ref{thm1.1} from Theorems~\ref{thm1.2},
\ref{thm1.3}, and~\ref{thm1.4}.

\begin{proof}[Proof of Theorem~\ref{thm1.1}]\leavevmode
  \begin{enumerate}[wide]
    \item Let $a,b\in\bar{\Q}$ be 
such that $a^{d}\neq b^{d}$ and $|S_{a,b}|>0$. By Theorems~\ref{thm1.2},
\ref{thm1.3}, and~\ref{thm1.4}, there exist real numbers
$c_{1},c_{2},c_{3}(\varepsilon),c_{4},\delta>0$
such that
\[
\max\{c_{1}h(a,b)-c_{2},\delta\}\leq\left\langle \boldsymbol{\mu}_{a},\boldsymbol{\mu}_{b}\right\rangle \leq\left(\varepsilon+\frac{c_{3}(\varepsilon)}{|S_{a,b}|}\right)(h(a,b)+c_{4}).
\]
Therefore, we have
\[
\max\left\{ \frac{c_{1}h(a,b)-c_{2}}{h(a,b)+c_{4}},\frac{\delta}{h(a,b)+c_{4}}\right\} \leq\varepsilon+\frac{c_{3}(\varepsilon)}{|S_{a,b}|}.
\]
When
\[
h(a,b)=\frac{c_{2}+\delta}{c_{1}},
\]
we have
\[
\frac{c_{1}h(a,b)-c_{2}}{h(a,b)+c_{4}}=\frac{\delta}{h(a,b)+c_{4}}=\frac{c_{1}\delta}{c_{1}c_{4}+c_{2}+\delta}.
\]
Since
\[
\frac{c_{1}h(a,b)-c_{2}}{h(a,b)+c_{4}}\text{ is increasing}\quad\text{and}\quad\frac{\delta}{h(a,b)+c_{4}}\text{ is decreasing}
\]
as functions of $h(a,b)$, we have
\[
\max\left\{ \frac{c_{1}h(a,b)-c_{2}}{h(a,b)+c_{4}},\frac{\delta}{h(a,b)+c_{4}}\right\} \geq\frac{c_{1}\delta}{c_{1}c_{4}+c_{2}+\delta}.
\]
When $\varepsilon$ is small enough, we have
\[
|S_{a,b}|\leq\frac{c_{3}(\varepsilon)}{\frac{c_{1}\delta}{c_{1}c_{4}+c_{2}+\delta}-\varepsilon}.
\]

\item Let $a\in\C\backslash\bar{\Q}$ and $b\in\C$ be 
such that $a^{d}\neq b^{d}$ and $|S_{a,b}|>0$. Since each $t\in S_{a,b}$
satisfies $f_{t}^{m}(a)=f_{t}^{n}(a)$ and $f_{t}^{k}(b)=f_{t}^{l}(b)$
for some $m>n\geq0$ and $k>l\geq0$, the field $\Q(a,b,S_{a,b})$
has transcendence degree $1$ over~$\Q$. We may view $\Q(a,b,S_{a,b})$
as the function field $K(X)$ of an algebraic curve $X$ defined over
a number field~$K$. For all but finitely many $x\in X(\bar{\Q})$,
the specializations satisfy $a(x)^{d}\neq b(x)^{d}$ and $t_{1}(x)\neq t_{2}(x)$
for any $t_{1},t_{2}\in S_{a,b}$. Therefore, the uniform bound for
the algebraic case also works for the complex case.
\hfill\qedhere
\end{enumerate}
\end{proof}



\begin{thebibliography}{DKY22+++}
  
\bibitem[BD11]{BD1}
M.~Baker and L.~DeMarco,
\emph{Preperiodic points and unlikely intersections},
Duke Math.~J.\ \textbf{159} (2011), no.~1, 1--29,
\doi{10.1215/00127094-1384773}.

\bibitem[BD13]{BD2}
\bysame, 
\emph{Special curves and postcritically finite polynomials},
Forum Math.\ Pi \textbf{1} (2013), Paper No.~e3,
\doi{10.1017/fmp.2013.2}.

\bibitem[BR06]{BR1}
M.~Baker and R.~Rumely,
\emph{Equidistribution of small points, rational dynamics, and potential theory}, Ann.\ Inst.\ Fourier (Grenoble) \textbf{56} (2006), no.~3, 625--688,
\doi{10.5802/aif.2196}.

\bibitem[BR10]{BR2}
\bysame,
\emph{Potential theory and dynamics on the Berkovich projective line},
Math.\ Surveys  Monogr., vol.~159,
Amer.\ Math.\ Soc., Providence, RI, 2010, \doi{10.1090/surv/159}.

\bibitem[Bea90]{Be}
A.\,F.~Beardon,
\emph{Symmetries of Julia sets},
Bull.\ London Math.\ Soc.\ \textbf{22} (1990), no.~6, 576--582,
\doi{10.1112/blms/22.6.576}.

\bibitem[Bil97]{Bi}
Y.~Bilu,
\emph{Limit distribution of small points on algebraic tori},
Duke Math.~J.\ \textbf{89} (1997), no.~3, 465--476,
\doi{10.1215/S0012-7094-97-08921-3}.

\bibitem[BFT18]{BFT}
F.~Bogomolov, H.~Fu, and Y.~Tschinkel,
\emph{Torsion of elliptic curves and unlikely intersections}, in:
\emph{Geometry and physics. Vol.~I}, pp.~19--37, Oxford Univ.\ Press, Oxford, 2018, \doi{10.1093/oso/9780198802013.003.0002}.

\bibitem[BH88]{BH}
B.~Branner and J.\,H.~Hubbard,
\emph{The iteration of cubic polynomials. I. The global topology of parameter space},
Acta Math.\ \textbf{160} (1988), no.~3-4, 143--206,
\doi{10.1007/BF02392275}.

\bibitem[CL06]{CL1}
A.~Chambert-Loir,
\emph{Mesures et \'equidistribution sur les espaces de Berkovich},
J.~reine angew.\ Math.\ \textbf{595} (2006), 215--235,
\doi{10.1515/crelle.2006.049}.

\bibitem[CL11]{CL2}
\bysame,
\emph{Heights and measures on analytic spaces. A survey of recent results, and some remarks},
in: \emph{Motivic integration and its interactions with model theory and non-Archimedean geometry. Vol.~II}, pp.~1--50,
London Math.\ Soc.\ Lecture Note Ser., vol.~384,
Cambridge Univ.\ Press, Cambridge, 2011,
\doi{10.1017/CBO9780511984433.002}.

\bibitem[DKY20]{DKY1}
L.~DeMarco, H.~Krieger, and H.~Ye,
\emph{Uniform Manin-Mumford for a family of genus 2 curves},
Ann.\ of Math.~(2) \textbf{191} (2020), no.~3, 949--1001,
\doi{10.4007/annals.2020.191.3.5}. 

\bibitem[DKY22]{DKY2}
\bysame, 
\emph{Common preperiodic points for quadratic polynomials},
J.~Mod.\ Dyn.\ \textbf{18} (2022), 363--413,
\doi{10.3934/jmd.2022012}.

\bibitem[DF08]{DF}
R.~Dujardin and C.~Favre,
\emph{Distribution of rational maps with a preperiodic critical point},
Amer.~J.\ Math.\ \textbf{130} (2008), no.~4, 979--1032,
\doi{10.1353/ajm.0.0009}.

\bibitem[FRL06]{FRL}
C.~Favre and J.~Rivera-Letelier,
\emph{\'Equidistribution quantitative des points de petite hauteur sur la droite projective}, 
Math.\ Ann.\ \textbf{335} (2006), no.~2, 311--361,
\doi{10.1007/s00208-006-0751-x}.

\bibitem[Fil17]{Fi}
P.~Fili,
\emph{A metric of mutual energy and unlikely intersections for dynamical systems}, preprint \arXiv{1708.08403v1} (2017). 

\bibitem[FP19]{FP}
P.~Fili and L.~Pottmeyer,
\emph{Quantitative height bounds under splitting conditions},
Trans.\ Amer.\ Math.\ Soc.\ \textbf{372} (2019), no.~7, 4605--4626,
\doi{10.1090/tran/7656}.

\bibitem[GHT13]{GHT1}
D.~Ghioca, L.-C.~Hsia, and T.~Tucker,
\emph{Preperiodic points for families of polynomials},
Algebra Number Theory \textbf{7} (2013), no.~3, 701--732,
\doi{10.2140/ant.2013.7.701}.

\bibitem[GHT15]{GHT2}
\bysame, 
\emph{Preperiodic points for families of rational maps},
Proc.\ Lond.\ Math.\ Soc.~(3) \textbf{110} (2015), no.~2, 395--427,
\doi{10.1112/plms/pdu051}.

\bibitem[MZ08]{MZ1}
D.~Masser and U.~Zannier,
\emph{Torsion anomalous points and families of elliptic curves},
C.\ R.\ Math.\ Acad.\ Sci.\ Paris \textbf{346} (2008), no.~9-10, 491--494,
\doi{10.1016/j.crma.2008.03.024}.

\bibitem[MZ10]{MZ2}
\bysame, 
\emph{Torsion anomalous points and families of elliptic curves},
Amer.~J.\ Math.\ \textbf{132} (2010), no.~6, 1677--1691,
\doi{10.1353/ajm.2010.a404145}.

\bibitem[MZ12]{MZ3}
\bysame,
\emph{Torsion points on families of squares of elliptic curves},
Math.\ Ann.\ \textbf{352} (2012), no.~2, 453--484,
\doi{10.1007/s00208-011-0645-4}.

\bibitem[PST12]{PST}
C.~Petsche, L.~Szpiro, and T.~Tucker,
\emph{A dynamical pairing between two rational maps},
Trans.\ Amer.\ Math.\ Soc.\ \textbf{364} (2012), no.~4, 1687--1710,
\doi{10.1090/S0002-9947-2011-05350-X}.

\bibitem[YZ13]{YZ}
X.~Yuan and S.~Zhang,
\emph{The arithmetic Hodge index theorem for adelic line bundles II},
preprint \arXiv{1304.3539v1} (2013). 

\bibitem[Zan12]{Za}
U.~Zannier,
\emph{Some problems of unlikely intersections in arithmetic and geometry} (with appendixes by D.~Masser),
Ann.\ of Math.\ Stud., vol.~181,
Princeton Univ.\ Press, Princeton, NJ, 2012,
\doi{10.23943/princeton/9780691153704.001.0001}.

\bibitem[Zha95]{Zh}
S.~Zhang,
\emph{Small points and adelic metrics},
J.~Algebraic Geom.\ \textbf{4} (1995), no.~2, 281--300.
\end{thebibliography}
\end{document}